\documentclass{gtart_h}  

\def\ifplaintex{\expandafter\ifx\csname documentclass\endcsname\relax}

\def\ifplaintex{\expandafter\ifx\csname documentclass\endcsname\relax}

\ifplaintex 
\else
\fi

\expandafter\ifx\csname epsfbox\endcsname\relax\input epsf\fi

\def\gt{{\mathsurround=0pt\it $\cal G\mskip-2mu$eometry \&\ 
$\cal T\!\!$opology}}        

\def\gtp{{\mathsurround=0pt\it $\cal G\mskip-2mu$eometry \&\ 
$\cal T\!\!$opology $\cal P\!$ublications}}  

\def\lognumber#1{\def\thelognumber{#1}}
\def\volumenumber#1{\def\thevolumenumber{#1}}
\def\papernumber#1{\def\thepapernumber{#1}}
\def\volumeyear#1{\def\thevolumeyear{#1}}

\def\pagenumbers#1#2{\def\startpage{#1}\def\finishpage{#2}}
\def\published#1{\def\publishdate{#1}}
\def\proposed#1{\def\theproposer{#1}}
\def\seconded#1{\def\theseconders{#1}}
\def\received#1{\def\receiveddate{#1}}
\def\revised#1{\def\reviseddate{#1}}
\def\accepted#1{\def\accepteddate{#1}}

\def\asciiaddress#1{\def\theasciiaddress{#1}}
\def\asciiemail#1{\def\theasciiemail{#1}}

\long\def\asciiabstract#1{\long\def\theasciiabstract{#1}}
\def\asciikeywords#1{\def\theasciikeywords{#1}}

\let\\\par\let\thelognumber\relax
\let\thevolumenumber\relax\let\thepapernumber\relax
\let\thevolumeyear\relax\let\thesamplenumber\relax\let\startpage\relax
\let\finishpage\relax\let\publishdate\relax\let\receiveddate\relax
\let\reviseddate\relax\let\accepteddate\relax\let\theasciititle\relax
\let\theasciiauthors\relax\let\theasciiaddress\relax
\let\theasciiabstract\relax\let\theasciikeywords\relax
\let\theasciiemail\relax\let\theshortauthors\relax\let\theshorttitle\relax

\long\def\maketitlep{   

\count0=\startpage

\gt\hfill      
\hbox to 77pt{\vbox to 0pt{\vglue -15pt\epsfbox{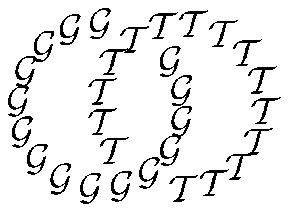}\vss}\hss}
\break
{\small\ifx\thesamplenumber\relax 
Volume \else Sample
\fi\thevolumenumber\ (\thevolumeyear)
\startpage--\finishpage\nl
Published: \publishdate}
\vglue 0.5truein plus 0.4fil minus 0.1truein

{\parskip=0pt\leftskip 0pt plus 1fil\def\\{\par\smallskip}{\ifplaintex\large
\else\Large\fi\bf\thetitle}\par\medskip}   

\vglue 0pt plus 0.1fil 

{\parskip=0pt\leftskip 0pt plus 1fil\def\\{\par}{\sc\theauthors}
\par\medskip}

\vglue 0pt plus 0.1fil 

{\small\parskip=0pt\let\newline\\
{\leftskip 0pt plus 1fil\def\\{\par}{\sl\theaddress}\par}
\expandafter\ifx\theemail\relax    
\relax\else\vglue 5pt plus 0.02fil minus 2pt\def\\{\stdspace{\rm 
and}\stdspace} 
\cl{Email:\stdspace\tt\theemail}\fi
\ifx\theurl\relax                  
\relax\else\vglue 5pt plus 0.02fil minus 2pt\def\\{\stdspace{\rm 
and}\stdspace}
\cl{URL:\stdspace\tt\theurl}\fi\par}

\vglue 7pt plus 0.3fil minus 3pt

{\bf Abstract}
\vglue 5pt plus 0.1fil minus 2pt

\theabstract

\vglue 7pt plus 0.3fil minus 3pt

{\bf AMS Classification numbers}\quad Primary:\quad \theprimaryclass

Secondary:\quad \thesecondaryclass

\vglue 5pt plus 0.3fil minus 2pt

{\bf Keywords:}\quad \thekeywords

\vglue 10pt plus 0.5fil minus 5pt

{\small  Proposed: \theproposer\hfill Received: \receiveddate\nl
Seconded: \theseconders\hfill 
\ifx\reviseddate\relax                         
Accepted: \accepteddate                        
\else
Revised: \reviseddate                          
\fi}
\eject
}       

\font\phead=cmsl9 scaled 950
\font=cmsl9 scaled 1050
\font\pnum=cmbx10 scaled 913
\font\lnum=cmbx10 
\font\pfoot=cmsl9 scaled 950
\font=cmsl9 scaled 1050
\ifplaintex
\headline{\vbox to 0pt{\vskip -4.5mm\line{\small\phead\ifnum
\count0=\startpage ISSN 1364-0380 (on line)
1465-3060 (printed) \hfill {\pnum\folio}\else\ifodd\count0\def\\{ }%
\ifx\theshorttitle\relax\thetitle\else\theshorttitle\fi\hfill{\pnum\folio}
\else\def\\{ and }{\pnum\folio}\hfill\ifx\theshortauthors\relax\theauthors
\else\theshortauthors\fi\fi\fi}\vss}}
\footline{\vbox to 0pt{\vglue 0mm\line{\small\pfoot\ifnum\count0=\startpage
\copyright\ \gtp\hfill\else
\gt, Volume \thevolumenumber\ (\thevolumeyear)\hfill\fi}\vss
}}
\else
\makeatletter
\def\@oddhead{{\small\lhead\ifnum\count0=\startpage ISSN 1364-0380 (on line)
1465-3060 (printed) \hfill {\lnum\number\count0}\else\ifodd\count0
\def\\{ }\ifx\theshorttitle\relax \thetitle \else\theshorttitle\fi\hfill
{\lnum\number\count0}\else\def\\{ and }{\lnum\number\count0}
\hfill\ifx\theshortauthors\relax 
\theauthors\else\theshortauthors\fi\fi\fi}}\def\@evenhead{\@oddhead}
\def\@oddfoot{\small\lfoot\ifnum\count0=\startpage\copyright\ \gtp\hfill\else
\gt, Volume \thevolumenumber\ (\thevolumeyear)\hfill\fi}
\def\@evenfoot{\@oddfoot}
\makeatother
\fi

\newwrite\gtoutfile
\long\gdef\makeheadfile{  
{\def\\{, }\def\s{ }
\immediate\openout\gtoutfile head.xxx
\immediate\write\gtoutfile{Proxy-for: \ifx\theasciiauthors\relax
\theauthors\else\theasciiauthors\fi\s<\ifx\theasciiemail\relax\theemail\else\theasciiemail\fi>}
\immediate\write\gtoutfile{\noexpand\\}
\immediate\write\gtoutfile{Authors: \ifx\theasciiauthors\relax
\theauthors\else\theasciiauthors\fi}
{\def\\{ }\immediate\write\gtoutfile{Title: \ifx\theasciititle\relax
\thetitle\else\theasciititle\fi}}
\immediate\write\gtoutfile{Subj-class: GT or SG or MG etc}
\immediate\write\gtoutfile{MSC-class: \theprimaryclass\ifx\thesecondaryclass\relax\else, \thesecondaryclass\fi}
\immediate\write\gtoutfile{Journal-ref: Geom. Topol. \thevolumenumber
(\thevolumeyear) \startpage-\finishpage}
\immediate\write\gtoutfile{Comments: Published by Geometry and Topology at}
\immediate\write\gtoutfile{\s\s http://www.maths.warwick.ac.uk/gt/GTVol\thevolumenumber/paper\thepapernumber.abs.html}
\immediate\write\gtoutfile{\noexpand\\}
\immediate\write\gtoutfile{}
\ifx\theasciiabstract\relax
\immediate\write\gtoutfile{\theabstract}\else
\immediate\write\gtoutfile{\theasciiabstract}\fi
\immediate\write\gtoutfile{}
\immediate\write\gtoutfile{\noexpand\\}
\immediate\write\gtoutfile{}
\immediate\closeout\gtoutfile}}  

\def\maketitlepage{\maketitlep\makeheadfile}
\let\maketitle\maketitlepage

\lognumber{477}
\received{15 July 2004}
\volumenumber{9}\papernumber{30}\volumeyear{2005}
\pagenumbers{1295}{1336}   
\revised{20 July 2005}
\published{26 July 2005}
\accepted{4 July 2005}
\proposed{Martin Bridson}
\seconded{Ralph Cohen, Thomas Goodwillie}

\usepackage{amsmath,amssymb,epsf,labelfig} 
\usepackage[all]{xy} 
 
\def\figref#1{\hyperlink{#1anchor}{Figure~\ref*{#1}}}
\def\anchor#1{\noindent\hypertarget{#1anchor}{\smash{$\phantom{99}$}}\newline}

\newtheorem{thm}{Theorem}[section] 
\newtheorem{lem}[thm]{Lemma} 
\newtheorem{prop}[thm]{Proposition} 
\newtheorem{cor}[thm]{Corollary} 
\newtheorem{Th}{Theorem}

\theoremstyle{remark}

\newcommand{\Dif}{\textrm{Diff}}  
\newcommand{\Aut}{\operatorname{Aut}}
\newcommand{\Out}{\operatorname{Out}}
\newcommand{\conn}{\operatorname{conn}}

\let\col\colon

\newcommand{\hy}{\kern-.1pt\vrule height2.8ptwidth4ptdepth-2.5pt\kern1pt}

\newcommand{\D}{\mathcal{D}}
\newcommand{\Sph}{\mathcal{S}}
\newcommand{\DS}{\mathcal{DS}}
\newcommand{\Z}{\mathbb{Z}}

\newcommand{\del}{\partial} 
\newcommand{\bdy}{\partial}
\newcommand{\cdotss}{\mathinner{\mkern1mu\cdotp\mkern-2mu\cdotp\mkern-2mu\cdotp}}
\newcommand{\tild}{\widetilde}

\newcommand{\semi}{\ltimes}
\newcommand{\sqdot}{\vrule height4.1pt width2.2pt depth -1.9pt \hskip5pt \ignorespaces}
\newcommand{\lowsqdot}{\vrule height3.6pt width2.2pt depth -1.4pt \hskip5pt \ignorespaces}
\newcommand{\smalltensor}{\otimes}
\newcommand{\emp}{\varnothing}

\newcommand{\rar}{\rightarrow} 
\newcommand{\inc}{\hookrightarrow}
\newcommand{\Rar}{\longrightarrow}

\begin{document}

\title{Stabilization for the automorphisms\\of free groups with boundaries}

\authors{Allen Hatcher\\Nathalie Wahl} 
\address{Mathematics Department, Cornell University, Ithaca NY 14853, USA}
\secondaddress{Aarhus University, 116 Ny Munkegade, 8000 Aarhus C, Denmark}
\asciiaddress{Mathematics Department, Cornell University, Ithaca NY 14853, 
USA\\and\\Aarhus University, 116 Ny Munkegade, 8000 Aarhus C, Denmark}
\gtemail{\mailto{hatcher@math.cornell.edu}{\rm\qua 
and\qua}\mailto{wahl@imf.au.dk}}
\asciiemail{hatcher@math.cornell.edu\\wahl@imf.au.dk}

\begin{abstract} 
The homology groups of the automorphism group of a free group are
known to stabilize as the number of generators of the free group goes
to infinity, and this paper relativizes this result to a family of
groups that can be defined in terms of homotopy equivalences of a
graph fixing a subgraph. This is needed for the second author's recent
work on the relationship between the infinite loop structures on the
classifying spaces of mapping class groups of surfaces and
automorphism groups of free groups, after stabilization and
plus-construction. We show more generally that the homology groups of
mapping class groups of most compact orientable 3-manifolds, modulo
twists along $2${\hy}spheres, stabilize under iterated connected sum with
the product of a circle and a $2${\hy}sphere, and the stable groups are
invariant under connected sum with a solid torus or a ball. These
results are proved using complexes of disks and spheres in reducible
$3${\hy}manifolds.
\end{abstract}

\asciiabstract{%
The homology groups of the automorphism group of a free group are
known to stabilize as the number of generators of the free group goes
to infinity, and this paper relativizes this result to a family of
groups that can be defined in terms of homotopy equivalences of a
graph fixing a subgraph. This is needed for the second author's recent
work on the relationship between the infinite loop structures on the
classifying spaces of mapping class groups of surfaces and
automorphism groups of free groups, after stabilization and
plus-construction. We show more generally that the homology groups of
mapping class groups of most compact orientable 3-manifolds, modulo
twists along 2-spheres, stabilize under iterated connected sum with
the product of a circle and a 2-sphere, and the stable groups are
invariant under connected sum with a solid torus or a ball. These
results are proved using complexes of disks and spheres in reducible
3-manifolds.}

\asciikeywords{Automorphism groups of free groups, homological stability, mapping class groups of 3-manifolds}

\maketitlepage

\section*{Introduction}\addcontentsline{toc}{section}{Introduction} 
There is by now a significant literature on homological stability
properties of various families of discrete groups, including matrix
groups, mapping class groups of surfaces, and automorphism groups of
free groups. The present paper is yet another contribution in this
direction, for a family of groups $ A^s_{n,k} $ that can be thought of
as relative versions of automorphism groups of free groups.  The
groups $A^1_{n,k}$ arose first in \cite{W} in studying the relation
between mapping class groups of surfaces and automorphism groups of
free groups, and applications to this work provide our main
motivation. The extra parameter $s$ in the family, necessary for the
proof of the stability theorem, has also interests of its own.  In the
present paper we exploit the fact that $ A^s_{n,k} $ can be realized
as a quotient of the mapping class group of a certain
$3${\hy}manifold, and our proof of stability applies to give a
homological stability result for the analogous quotients of the
mapping class groups of general compact orientable $3${\hy}manifolds.

\begin{figure}[ht!]\small\anchor{fig1}
\SetLabels
(.15*.72) $n$\\
(.85*.68) $k$\\
\endSetLabels
\centerline{
\AffixLabels{\epsfbox{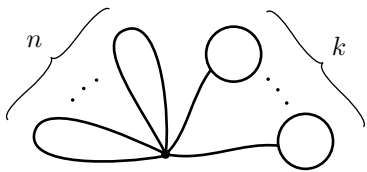}}} 
\caption{The graph $G_{n,k}$}\label{fig1}
\end{figure}

The quickest definition of $ A_{n,k}=A_{n,k}^1$ is in terms of
homotopy equivalences of graphs. Let $G_{n,k}$ denote the graph shown
in \figref{fig1} consisting of a wedge of $n$ circles together with
$k$ distinguished circles joined by arcs to the basepoint.

Then $A_{n,k}$ is the group of path-components of the space of
homotopy equivalences of $G_{n,k}$ that preserve the basepoint and
restrict to the identity on each of the $k$ distinguished circles. If
this last condition is relaxed to allow homotopy equivalences that
rotate each of the $ k $ distinguished circles, the group $ A_{n,k} $
is replaced by the group $ \Aut_{n,k} $ of automorphisms of the free
group on $ n + k $ generators that take each of the last $ k $
generators to a conjugate of itself. By restricting the more general
homotopy equivalences of $ G_{n,k} $ to their rotations of the
distinguished circles, we obtain a fibration whose associated long
exact sequence of homotopy groups is a short exact sequence
$$
1 \Rar \Z^k \Rar A_{n,k} \Rar \Aut_{n,k} \Rar 1
$$
expressing $ A_{n,k} $ as a central extension of $ \Aut_{n,k} $.  The
groups $ \Aut_{n,k} $ can be regarded as interpolations between
automorphism groups of free groups, which are the groups $ \Aut_{n,0}
$, and the groups $ \Aut_{0,k} $ which are the `symmetric'
automorphism groups studied for example in \cite{C} and \cite{Mc}. A
purely algebraic description of $ A_{n,k} $ is given in \cite{JW}.

The group $A_{n,k}$ relates to mapping class groups of surfaces in the
following way. Let $S$ be a surface of genus $g$ with $k+1$ boundary
components and let $\Gamma_{g,k+1}$ be its mapping class group, the
group of path components of the space of diffeomorphisms of $S$ fixing
the boundary pointwise.  There is a natural embedding of the graph
$G_{2g,k}$ in $S$ with the basepoint of $G_{2g,k}$ on the first
boundary component of $S$ and the $k$ distinguished circles of $G$
identified with the other $k$ boundary components of $S$, such that
the embedding is a homotopy equivalence relative to the basepoint and
the $k$ distinguished circles. This defines a map from
$\Gamma_{g,k+1}$ to $A_{2g,k}$, which is actually injective.  Dehn
twists along the last $k$ boundary circles of $S$ are identified under
this inclusion with the copy of $\Z^k$ in $A_{2g,k}$ arising in the
short exact sequence above.

The natural inclusions $ G_{n,k} \inc G_{n+1,k} $ and $ G_{n,k} \inc
G_{n,k+1} $ obtained by attaching an extra circle or circle-plus-arc
induce stabilization maps $ A_{n,k} \rar A_{n+1,k} $ and $ {A_{n,k}
\rar A_{n,k+1}} $ by extending homotopy equivalences via the identity
on the added parts. These stabilization maps are easily seen to be
injective. The main goal of the paper is to prove the following
result:

\begin{Th}\label{thmA}
The stabilization maps $ A_{n,k} \rar A_{n+1,k} $ and $ {A_{n,k} \rar
A_{n,k+1}} $ induce isomorphisms $ H_i(A_{n,k}) \rar H_i(A_{n+1,k}) $
and $ H_i(A_{n,k}) \rar H_i(A_{n,k+1}) $ for $n \ge 3i+3 $ and $k \ge
0$.
\end{Th}

In particular, the case $k=0$ in the first stabilization gives a new
proof of homology stability for automorphism groups of free groups,
although the ranges $n\ge 2i+3$ obtained in \cite{HV1} and $ n \ge
2i+2$ in \cite{HV3} are better than the one we obtain here.  The proof
of the theorem also shows that the quotient map $\Aut(F_n) \rar
\Out(F_n)$ induces an isomorphism on $ H_i $ for $ n \ge 3i+3
$. Again, a better range is given in \cite{HV3}.

The theorem says that the stable $ A_{n,k} $ groups have the same
homology as the stable automorphism groups of free groups. Not a great
deal is known about these stable homology groups. In low dimensions
they are trivial rationally \cite{HV1}, \cite{HV2}, and this is
conjectured to hold in all dimensions. The only nontrivial torsion
known so far is a copy of the homology of the infinite symmetric
group, or equivalently $ \Omega^{\infty}S^{\infty} $, as noted in
\cite{H3}.

The above theorem is one of the main ingredients in \cite{W} in
showing that the natural map from the stable mapping class group of
surfaces to the stable automorphism group of free groups induces an
infinite loop map on the classifying spaces of these groups after
plus-construction.  The infinite loop space structure of the stable
mapping class group is governed by a cobordism category, the surfaces
being cobordisms between disjoint unions of circles \cite{T}.  Graphs
with distinguished circles, thought of as {\em boundary circles}, can
be used to define a generalized cobordism category where graphs as
well as surfaces define cobordisms and where gluing along boundary
circles is compatible in an appropriate sense with the map
$\Gamma_{g,k+1}\rar A_{2g,k}$ described above. This category is used
to compare the two infinite loop spaces and homological stability for
the groups $A_{n,k}$ is needed for identifying the homotopy type of
this generalized cobordism category.  From \cite{MW} and \cite{Gal}
the homology of the stable mapping class group is known (and it is
very rich!). Little is known however about the map from the mapping
class groups to the automorphism groups. If the map could be shown to
be rationally trivial, for example, this would give further evidence
for the conjecture that the stable rational homology of the
automorphism groups is trivial.

To prove the theorem we actually need to work with the more general
groups $A_{n,k}^s$ which can be defined as follows.  Consider the
graph $G_{n,k}^s$ obtained from $G_{n,k}$ by wedging $s-1$ edges at
the basepoint.  We think of the free ends of the new edges as extra
basepoints. (When $s=0$, we forget the basepoint of $G_{n,k}$.) Define
$A_{n,k}^s$ as the group of path components of the space of homotopy
equivalences of $G_{n,k}^s$ which fix the $k$ distinguished circles as
well as the $s$ basepoints.  There is a short exact sequence
$$
1\Rar F_{n+k}\Rar A_{n,k}^s\Rar A_{n,k}^{s-1} \Rar 1
$$
where the map $A_{n,k}^s\rar A_{n,k}^{s-1} $ forgets the last
basepoint. When $s\ge 2$, this map splits and in fact
$A_{n,k}^s=(F_{n+k})^{s-1}\semi A_{n,k}^1$ where $A_{n,k}^1$ acts
diagonally on $(F_{n+k})^{s-1}$ via the map $A^1_{n,k}\rar
\Aut(F_{n+k})$.  Particular cases of these groups are
$A_{n,0}^0=\Out(F_n)$, $A_{n,0}^1=\Aut(F_n)$ and $A_{n,0}^2=F_n\semi
\Aut(F_n)$.

We show in Theorem~\ref{proc4.1} that the stabilization map $
{A^s_{n,k} \rar A^{s+1}_{n,k}} $ and the map $ {A^s_{n,k} \rar
A^{s-1}_{n,k}} $ for $s\ge 1$ also induce isomorphisms in homology
when $n\ge 3i+3$.  A corollary of the stability for $A_{n,0}^2\rar
A_{n,0}^1$ is that the twisted homology group $H_i(\Aut(F_n),{\bf
Z}^n)$ is trivial when $n\ge 3i+9$. This is immediate from the
Leray--Hochschild--Serre spectral sequence for the short exact
sequence displayed above, with $ s = 2 $ and $ k = 0 $.

Our proof of homological stability for $A_{n,k}$, and in fact
$A_{n,k}^s$, is based on an interpretation of these groups in terms of
the mapping class group of compact $3${\hy}manifolds $M^s_{n,k}$
obtained from the connected sum of $n$ copies of $S^1 \times S^2$ by
deleting the interiors of $s$ disjoint balls and $k$ disjoint
unknotted, unlinked solid tori $S^1 \times D^2$. Equivalently, $
M^s_{n,k}$ is the connected sum of $n$ copies of $ S^1 \times S^2 $,
$s$~balls, and $k$ solid tori.  If $\Dif(M^s_{n,k})$ denotes the group
of orientation-preserving diffeomorphisms of $M^s_{n,k}$ that restrict
to the identity on the boundary, then $A^s_{n,k}$ is the quotient of
the mapping class group $ \pi_0\Dif(M^s_{n,k})$ by the subgroup
generated by Dehn twists along embedded spheres.  This subgroup is
rather small, the product of finitely many $\Z_2$'s.

It seemed worthwhile to write the proof in its most natural level of
generality. Thus we actually consider manifolds $ M^s_{n,k} $ obtained
from a fixed compact orientable manifold $ N $ by taking the connected
sum with $n$ copies of $S^1 \times S^2$, $s$~balls, and $ k $ solid
tori, where we assume that none of the connected summands of $N$ has
universal cover a counterexample to the Poincar\'e conjecture. We let
now $ A^s_{n,k} $ denote the mapping class group $
\pi_0\Dif(M^s_{n,k}) $ with twists along $2${\hy}spheres factored
out. Again these generate just a product of finitely many $ \Z_2
$'s. When $ s \ge 1 $ there are natural stabilizations with respect to
$ n $, $ k $, and $ s $ by enlarging $ M^s_{n,k} $ and extending
diffeomorphisms by the identity. The more general form of the
preceding theorem is then:

\begin{Th}\label{thmB}
The maps $ H_i(A^s_{n,k}) \rar H_i(A^s_{n+1,k}) $, $ H_i(A^s_{n,k})
\rar H_i(A^s_{n,k+1}) $, and $ H_i(A^s_{n,k}) \rar H_i(A^{s+1}_{n,k})
$ induced by stabilization are isomorphisms for $n \ge 3i+3 $, $k \ge
0 $, and $ s \ge 1 $.\end{Th}

One might hope that these stability results could be part of a program
to extend the work of Madsen--Weiss on mapping class groups of
surfaces to $3${\hy}manifolds, but there are some obstacles to doing
this. For surfaces the contractibility of the components of $ \Dif $
was an essential step in the Madsen--Weiss program, allowing $
\pi_0\Dif $ to be replaced by $ \Dif $. However, the components of $
\Dif(M^s_{n,k}) $ are far from contractible, according to the results
of \cite{HL}, and the homology of $ \Dif(M^s_{n,k}) $ itself does not
stabilize. This suggests that for $3${\hy}manifolds a different sort
of stabilization than by connected sum might be better if one wanted
to use the Madsen--Weiss techniques, a stabilization that stays within
the realm of irreducible $3${\hy}manifolds where $ \Dif $ often has
contractible components, as shown in \cite{H1} and \cite{Gab}.

Here is how the paper is organized. The first section is devoted to
showing the definition of $A^s_{n,k}$ in terms of mapping class groups
agrees with the graph-theoretic definition. The second section proves
some preliminary properties of systems of spheres and disks in a
$3${\hy}manifold, generalizing \cite{H3}. These properties are used in
the third section to prove that the complexes of disks and spheres we
are interested in are highly connected. Then in the fourth section the
spectral sequence arguments proving homological stability are given.

We have tried to make the paper largely independent of \cite{H3}, but
the reader might still find it helpful to consult this earlier paper
at certain points.

The second author was supported by a Marie Curie Fellowship of the
European Community under contract number HPMF-CT-2002-01925.

\section{Diffeomorphism groups} 

In the introduction we first defined $A^s_{n,k}$ as the group of
path-components of the space of homotopy equivalences of the graph
$G^s_{n,k}$ that fix the $s$ basepoints and the $k$ distinguished
circles. A homotopy equivalence satisfying this condition is in fact a
homotopy equivalence relative to the basepoints and distinguished
circles, according to Proposition~0.19 of \cite{H4}, so $A^s_{n,k}$ is
indeed a group, with inverses. When $s\ge 1$, by assigning to each
homotopy equivalence its induced automorphism of $ \pi_1 G^s_{n,k} =
F_{n+k} $ one obtains a short exact sequence
$$
1 \Rar \Z^k\times (F_{n+k})^{s-1} \Rar A^s_{n,k} \Rar \Aut_{n,k} \Rar
1
$$
The copy of $\Z^k$ in the kernel is generated by the homotopy
equivalences produced by rotating one of the distinguished circles
through 360 degrees and dragging along the arc connecting this circle
to the basepoint so that this arc wraps once around the circle. The
$i$th copy of $F_{n+k}= \pi_1G_{n,k}^s$ is identified with the
homotopy classes of maps from the interval to $G_{n,k}^s$ mapping the
endpoints to the endpoints of the $i$th extra edge.

Our aim in this section is to relate $A^s_{n,k}$ to the mapping class
group of $ M^s_{n,k} $, where $M^s_{n,k}$ is the connected sum
$$
M^s_{n,k} = \bigl(\#_n(S^1\times S^2)\bigr)\ \#\ \bigl(\#_k(S^1\times
D^2)\bigr)\ \#\ \bigl(\#_sD^3\bigr)
$$
For an orientable $ 3 ${\hy}manifold $ M $ we denote by $ \Dif(M) $
the group of orientation-preserving diffeomorphisms of $ M $ that fix
its boundary pointwise. When the boundary of $ M $ is nonempty, as
will usually be the case in this paper, diffeomorphisms that restrict
to the identity on the boundary automatically preserve orientation.

\begin{thm}\label{proc1.1}
There is an exact sequence
$$
1 \Rar K^s_{n,k} \Rar \pi_0\Dif(M^s_{n,k}) \Rar A^s_{n,k} \Rar 1
$$ 
where the kernel $K^s_{n,k}$ is the subgroup of $ \pi_0\Dif(M^s_{n,k})
$ generated by Dehn twists along embedded $2${\hy}spheres. This
subgroup is a product of at most $n+k+s$ copies of $\Z_2$.
\end{thm}

In the case $k=0$ and $s=0,1$ when $ A^1_{n,0} = \Aut(F_n) $ and
$A_{n,0}^0= \Out(F_n)$, this is a theorem of Laudenbach
\cite[III.4.3]{L} who showed that $ K^0_{n,0}=K_{n,0}^1 $ is the
product of exactly $n$ copies of $\Z_2$. (Laudenbach considered
diffeomorphisms of $ M^0_{n,0}$ fixing a point rather than
diffeomorphisms of $M^1_{n,0}$ fixing the boundary, but the group of
isotopy classes is the same in both cases.) For our purposes we only
need to know that $K^s_{n,k}$ is generated by twists along
$2${\hy}spheres. As we will see, twists along spheres act trivially on
embedded spheres and disks, up to isotopy, so the natural action of
the diffeomorphism group of $ M $ on the complexes of isotopy classes
of spheres and disks that we define in Section 3 will induce an action
of $ A^s_{n,k} $ on these complexes.

\begin{proof}  We first consider the case when $s=1$. As before, we drop $s=1$ from the notation so that $A_{n,k}=A_{n,k}^1$ etc.
The first step will be to construct a homomorphism $ \alpha \col
\pi_0\Dif(M_{n,k}) \rar A_{n,k} $. To do this we first describe a
natural embedding $ G_{n,k} \inc M_{n,k} $. We can construct $ M_{n,k}
$ in the following way. Start with the space $X$ obtained from a ball
$ B $ by removing the interiors of $ 2n+k $ disjoint subballs $ B_i $,
$i = 1, \cdotss , 2n+k $. Then $M_{n,k}$ is the quotient space of $X$
obtained by identifying $ \bdy B_{2i-1} $ with $ \bdy B_{2i} $ for $ i
= 1,\cdotss ,n $ to form nonseparating spheres $ S_i \subset M_{n,k} $
and identifying two disjoint disks $ D'_i $ and $ D''_i $ in $ \bdy
B_i $ for $ {i = 2n+1, \cdotss ,2n+k} $ to form disks $D_i \subset
M_{n,k}$.
\begin{figure}[ht!]\anchor{fig2}
\centerline{\epsfbox{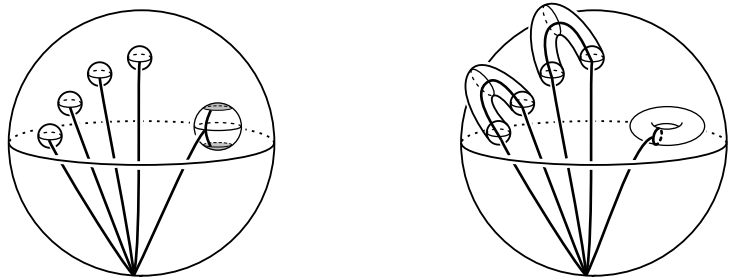}}
\caption{The graph $G_{n,k}$ in $M_{n,k}$}\label{fig2}
\end{figure}
Construct a tree $T \subset X$ by taking arcs $ a_i $ joining a
basepoint $x_0 \in \bdy B$ to points $ x_i \in \bdy B_i $ (with $x_i$
disjoint from the disks $D'_i$ and $D''_i$ for $ {i = 2n+1, \cdotss ,
2n+k} $), together with arcs $b_i$ in $ \bdy B_i $ for $ i = 2n+1,
\cdotss , 2n+k $ joining points $ x'_i \in \bdy D'_i$ and $x''_i \in
\bdy D''_i$ and passing through $x_i$.  We may assume the points
$x_i$, $x'_i$, and $x''_i$ are chosen so that the pairs $
(x_{2i-1},x_{2i}) $ for $ i = 1,\cdotss ,n $ and the pairs $
(x'_i,x''_i) $ for $ i = 2n+1, \cdotss , 2n+k $ match up under the
identifications which form $M_{n,k}$, and then the image of $T$ in
$M_{n,k}$ is a copy of $G_{n,k}$. Let $i \col G_{n,k} \rar M_{n,k} $
be this embedding.  With this embedding the $k$ distinguished circles
of $G_{n,k}$ lie in $ \bdy M_{n,k}$ so we call them {\em boundary
circles}. There is a retraction $ X \rar T $ obtained by first
deformation retracting $X$ onto the union of $T$ with the spheres
$\bdy B_i$ and then retracting this union onto $T$ by collapsing each
sphere $\bdy B_i$ to its intersection with $T$, which is either a
point $x_i$ or an arc $b_i$. These retractions can be taken to respect
the identifications which produce $M_{n,k}$, so they induce a
retraction $r \col M_{n,k} \rar i(G_{n,k}) $. This retraction is not a
homotopy equivalence but it does induce an isomorphism on $ \pi_1 $,
as does the inclusion map $i$.

Now we are able to define a map $ \alpha \col \pi_0\Dif(M_{n,k}) \rar
A_{n,k} $ by sending a diffeomorphism $ \varphi \col M_{n,k} \rar
M_{n,k} $ to the map $ r \varphi i \col G_{n,k} \rar G_{n,k} $. This
composition is a homotopy equivalence since it induces an isomorphism
on $\pi_1$. Since it fixes the basepoint and the boundary circles, it
therefore represents an element of~$A_{n,k}$.

To show that $\alpha$ is a homomorphism, we need the following fact:
For every map $j \col G_{n,k} \rar M_{n,k} $ that agrees with $i$ on
the basepoint and boundary circles there is a homotopy $ irj \simeq j
$ relative to the basepoint and boundary circles. If this is true then
by taking $ j = \varphi i $ we obtain a homotopy $ (r \psi i)(r
\varphi i) \simeq r(\psi\varphi) i $. Hence $\alpha$ is a
homomorphism.

To verify the fact, it suffices to homotope $j$ to have image in $
i(G_{n,k})$ by a homotopy relative to the basepoint and boundary
circles. Such a homotopy exists since $ i $ induces a surjection on $
\pi_1 $, so the pair $ \bigl(M_{n,k},i(G_{n,k})\bigr) $ is $ 1
${\hy}connected, which means that any map $ (I, \bdy I) \rar
\bigl(M_{n,k},i(G_{n,k})\bigr) $ is homotopic, fixing $ \bdy I $, to a
map into $ i(G_{n,k}) $.

The next step of the proof is to show that $\alpha$ is
surjective. Following the method of \cite{L} this will be done by
showing that generators for $ A_{n,k} $ can be realized by
diffeomorphisms of $M_{n,k}$. According to \cite{JW} the group
$A_{n,k}$ has the following generators, described by their effect on a
basis $\{x_1,\cdotss,x_n,y_1,\cdotss,y_k\}$ for $F_{n+k} =
\pi_1(G_{n,k})$, where we write the action on generators only when it
is nontrivial.

\smallskip
\halign{\quad #\hfil & \quad #\hfil & \quad #\hfil\cr $P_{i,j}$ & $x_i
\rar x_j $ and $ x_j \rar x_i $ & $ 1 \le i,j \le n $, $ i\ne j$ \cr
$I_i$ & $x_i \rar x^{-1}_i $ & $ 1 \le i \le n$ \cr $(x_i;x_j)$ & $x_i
\rar x_i x_j $ & $ 1 \le i,j \le n $, $ i\ne j$ \cr $(x_i;y_j)$ & $x_i
\rar x_i y_j $ & $ 1 \le i \le n $, $ 1 \le j \le k $ \cr
$(x^{-1}_i;y_j)$ & $x_i \rar y^{-1}_j x_i $ & $ 1 \le i \le n $, $ 1
\le j \le k $ \cr $(y^\pm_i;x_j)$ & $y_i \rar x^{-1}_j y_i x_j $ & $ 1
\le i \le k $, $ 1 \le j \le n $ \cr $(y^\pm_i;y_j)$ & $y_i \rar
y^{-1}_j y_i y_j $ & $ 1 \le i,j \le k $ \cr }

When $i=j$ in the last set of generators one has the trivial
automorphism of $F_{n+k}$ but a nontrivial element of $A_{n,k}$, a
generator of the subgroup $\Z^k$.

The first three types of generators involve only the $x_i$'s and can
be realized as in \cite{L}. Thinking of the $S^1 \times S^2$ summands
of $M_{n,k}$ as handles, the generator $P_{i,j}$ interchanges the
$i$th and $j$th handles, $I_i$ switches the two ends of the $i$th
handle, and $(x_i;x_j)$ slides one end of the $i$th handle over the
$j$th handle. The result of this slide is a Dehn twist along a torus
that encloses one end of the $i$th handle and passes over the $j$th
handle. The next two types of generators, $(x_i;y_j)$ and
$(x^{-1}_i;y_j)$, can be realized by similar slides of one end of the
$i$th handle along a loop that represents $y_j$, going through the
$j$th $S^1 \times D^2$ summand of $M_{n,k}$. For the last two types of
generators we think of the $i$th $S^1 \times D^2$ summand as a
punctured solid torus attached to the rest of $M_{n,k}$ along a sphere
that gives the connected sum operation, and then we slide this
punctured solid torus over the $j$th handle to realize the conjugation
$(y^\pm_i;x_j)$, or through the $j$th $S^1 \times D^2$ summand to
realize $(y^\pm_i;y_j)$ for $ i \ne j$. Like the slide producing
$(x_i;x_j)$, these slides also yield Dehn twists along tori. The last
remaining type of generator $(y^\pm_i;y_i)$ is realized by a Dehn
twist along a torus parallel to the $i$th boundary torus of $M_{n,k}$,
with the twisting in the direction of the $S^1$ factor of the summand
$S^1 \times D^2$.

The last step of the proof is to describe the kernel of $\alpha$,
using the following commutative diagram
\vskip-8pt
$$\SelectTips{cm}{} \xymatrix{1\ar[r]& K_{n,k} \ar[r] &
\pi_0\Dif(M_{n,k}) \ar[r]^-\alpha \ar[dr] & A_{n,k}\ar[r] \ar[d] & 1\\
& & & \Aut(F_{n+k}) & }$$
\vskip-4pt where the top row is exact and the diagonal map is the
action on $\pi_1$. Let $\varphi$ be a diffeomorphism representing a
class in the kernel $K_{n,k}$, so in particular $\varphi$ induces the
identity on $ \pi_1(M_{n,k})$. By \cite[ A-III.3.2]{L} we know that
$\varphi$ also induces the identity on $\pi_2(M_{n,k})$. Let $\Sigma$
be a system of $n+k$ spheres in $M_{n,k}$ consisting of the $n$
spheres $S_i$ in the $S^1 \times S^2$ summands together with spheres
$S_{n+1},\cdotss,S_{n+k}$ separating off the $S^1 \times D^2$
summands. Splitting $M_{n,k}$ along $\Sigma$ produces a ball with
$2n+k$ punctures, along with $k$ copies of a punctured solid
torus. Because $\varphi$ induces the identity on $\pi_2$, it can be
isotoped so that it is the identity on $\Sigma$, using Laudenbach's
homotopy-implies-isotopy theorem for $2${\hy}spheres, taking care of
the spheres $S_i$ one by one. (For the induction step one uses
\cite[V.4.2]{L} which says that homotopic spheres that are disjoint
from a sphere $S$ and not isotopic to $S$ are still isotopic in the
complement of $S$; this is reproved in our Theorem~\ref{proc2.2}.) So
we can consider $\varphi$ separately on the components of
$M-\Sigma$. On the punctured ball it is isotopic to a product of Dehn
twists along a subset of the $2n+k$ boundary spheres of the
punctures. According to Lemma~\ref{proc1.2} which follows below, the
mapping class group of a punctured solid torus is $\Z_2 \times \Z$
with the first factor generated by a Dehn twist along the boundary
sphere and the second factor generated by a Dehn twist along the
boundary torus, in the $S^1$ direction. The latter twist is the
diffeomorphism realizing one of the generators $(y_i;y_i)$ of the
$\Z^k$ subgroup of $A_{n,k}$, so only the $\Z_2$ factor lies in the
kernel $K_{n,k}$.

The conclusion of all this is that the diffeomorphism $\varphi$
representing an element of the kernel of $ \alpha $ is isotopic to a
product of Dehn twists along a subset of the spheres of $\Sigma$. The
proof of the theorem for $s=1$ is then completed by observing that all
twists along $2${\hy}spheres lie in the kernel of $ \alpha $ since for
a given $2${\hy}sphere, the map $ i\col G_{n,k} \rar M_{n,k} $ can be
homotoped to intersect this sphere only in the fixed points of a twist
along the sphere.

An argument for deducing the theorem for arbitrary $ s $ from the case
$ s = 1 $ is given in Proposition~1 of \cite{HV3} when $k=0$, and the
same argument works equally well for arbitrary $k$.
\end{proof}

\begin{lem}\label{proc1.2}
$\pi_0\Dif(M^1_{0,1}) \cong \Z_2 \times \Z$ where the $\Z_2$ is
generated by a Dehn twist along the boundary sphere and the $\Z$ is
generated by a Dehn twist along the boundary torus in the direction of
the $S^1$ factor of the summand $ S^1 \times D^2$ of $M^1_{0,1} = (S^1
\times D^2) \,\,\#\ D^3$.
\end{lem}

\begin{proof}  Consider the fibration
$$
\Dif(M^1_{0,1}) \Rar \Dif(S^1 \times D^2) \Rar E(D^3,S^1 \times D^2)
$$
where $E(D^3,S^1 \times D^2)$ is the space of embeddings of a ball in
the interior of $S^1 \times D^2$ and the projection to $E$ is obtained
by restricting diffeomorphisms to a chosen ball in the interior of
$S^1 \times D^2$. Since elements of $\Dif(S^1 \times D^2)$ fix the
boundary pointwise, they preserve orientation, so we may as well take
$E(D^3,S^1 \times D^2)$ to consist of orientation-preserving
embeddings.

The total space $\Dif(S^1 \times D^2)$ of the fibration is
contractible. This can be seen by looking at a second fibration
$$
\Dif(D^3) \Rar \Dif(S^1 \times D^2) \Rar E(D^2,S^1 \times D^2)
$$
whose projection map is restriction to a meridian disk. The base space
consists of embeddings of a disk in $S^1 \times D^2 $ fixed on the
boundary. It is a classical fact that this space is path-connected,
and in fact it is contractible by \cite{H1}. The fiber consists of
diffeomorphisms of $S^1 \times D^2$ that are fixed on the boundary and
a meridian disk, so this can be identified with $\Dif(D^3)$ which is
contractible by \cite{H2}. Hence $\Dif(S^1 \times D^2)$ is also
contractible.

The base space of the first fibration is homotopy equivalent to $S^1
\times SO(3)$. The long exact sequence of the fibration then gives an
isomorphism
$$
\pi_0\Dif(M^1_{0,1}) \cong \pi_1\bigl(S^1 \times SO(3)\bigr) \cong \Z
\times \Z_2
$$
where the $Z_2$ corresponds to a Dehn twist about the boundary sphere
of $M^1_{0,1}$ and the $\Z$ is generated by the diffeomorphism
resulting from an isotopy dragging a ball in the interior of $S^1
\times D^2$ around the $S^1$ factor. This diffeomorphism is equivalent
to a Dehn twist about the boundary torus of $M^1_{0,1}$. \end{proof}

\section{Normal form} 

We are primarily interested in the manifolds $M^s_{n,k}$, but it will
be convenient and not much extra work to consider a more general
compact manifold $M$ that is the connected sum of $ n $ copies of $
{S^1 \times S^2} $, $ s $ balls, and $ k $ irreducible orientable
manifolds $Q_i$ that are neither a ball nor a counterexample to the
Poincar\'e conjecture. The latter condition is needed since we will be
using Laudenbach's theorem that homotopy implies isotopy for embedded
spheres, which requires this hypothesis. We will in fact assume that
no $ Q_i $ has universal cover a counterexample to the Poincar\'e
conjecture, although it is possible that this stronger assumption
could be avoided by taking extra care at certain points in the proofs.

We will study $ \pi_0\Dif(M) $ by studying isotopy classes of embedded
spheres and disks in $ M $. By an embedded sphere or disk we mean a
submanifold of $ M $ that is diffeomorphic to a sphere or disk, but
without a specific diffeomorphism being chosen. An isotopy of an
embedded sphere or disk is a smooth one-parameter family of such
submanifolds. This is the same as choosing a diffeomorphism from the
standard sphere or disk to the given submanifold and varying this map
to $M$ by isotopy. More generally, if we allow the map to vary by
homotopy instead of isotopy we can speak of two embedded spheres or
disks being homotopic as well as isotopic.

Besides individual spheres and disks, we also consider systems of
finitely many embedded spheres and disks in $M$. The spheres and disks
in a system are assumed to be disjoint, except that the boundary
circles of different disks are allowed to coincide. We will assume the
disks have their boundaries contained in a fixed finite collection $ C
$ of disjoint circles in $ \bdy M $. In the end we will only need the
case that $ C $ is a single circle, either in a boundary sphere of $ M
$ or in the boundary torus of a solid torus summand $ Q_i = S^1 \times
D^2 $, with the circle being a meridian $\{x\} \times \bdy D^2$ in the
latter case. However, to prove the results in Section 3 we will have
to allow larger collections $ C $. It will suffice to assume:

\items
\item[\sqdot] $C$ consists of at most two circles in each sphere
component of $ \bdy M $ and at most one circle in the boundary torus
of each summand $ Q_i = S^1 \times D^2 $, this circle being a
meridian.  \enditems

The reason for allowing two circles in a boundary sphere is that we
will be cutting $ M $ open along systems of disks with boundaries in
$C$. Cutting $M$ along a disk bounded by a circle of $C$ that is a
meridian in a boundary torus replaces the torus by a boundary sphere
containing two copies of the original circle of $C$ (bounding two
copies of the disk). On the other hand, cutting $M$ along a disk
bounded by a circle of $C$ in a boundary sphere replaces the sphere by
two spheres, each containing a copy of the original circle of $C$, so
the number of circles in a boundary sphere does not increase.

The spheres and disks in a system will be assumed to be
nontrivial. For spheres this means that they do not bound balls in $M$
and are not isotopic to boundary spheres of $M$ disjoint from
$C$. Nontriviality for disks means they are not isotopic, fixing their
boundary, to disks in $ \bdy M $ containing no circles of $ C $ in
their interior. Another assumption we make is that no two spheres in a
system are isotopic, and no two disks are isotopic fixing their
boundaries. This is equivalent to saying that no two spheres bound a
product $ S^2 \times I $ in $ M $ and no two disks bound a ball.

To save words in the rest of the paper we will always assume that
isotopies and homotopies of disks in a system fix the boundaries of
these disks.

Our purpose in this section is to develop a notion of normal form for
systems of spheres and disks in the manifold $ M $, generalizing
\cite{H3} which dealt with systems of spheres in the manifolds $
M^s_{n,0} $. The simplest definition of normal form is to say that a
system $S$ of spheres and disks in $M$ is in {\em normal form\/} with
respect to a fixed maximal sphere system $\Sigma$ if $S$ is transverse
to $\Sigma$ and the number of circles of $S \cap \Sigma$ is minimal
among all systems isotopic to $S$ and transverse to $\Sigma$. As we
will see, in the case of the manifolds $ M^s_{n,0} $ this boils down
to almost the same thing as the more complicated definition given in
\cite{H3}, with only a minor difference in how spheres of $S$ that are
isotopic to spheres of $\Sigma$ are treated.

The advantage of the present definition is that the existence of a
system in normal form within each isotopy class of systems is
obvious. This allows us to focus on what is the real point of normal
form, which is to obtain a reasonable uniqueness statement. The
strongest uniqueness statement would be that if two normal form
systems are isotopic then they are isotopic through normal form
systems. However, this could hold only if one first eliminated certain
knotting and linking phenomena. Fortunately knotting and linking are
not really an issue in view of the fact that Laudenbach's
homotopy-implies-isotopy theorem for spheres and sphere systems holds
also for systems of spheres and disks, as we show in
Theorem~\ref{proc2.2} below. We can then avoid knotting and linking
problems by making the following definition. Two systems $S$ and $S'$
transverse to $\Sigma$ are called {\em equivalent} if there is a
homotopy $S_t$ from $S$ to $S'$ through immersions such that $S_t$ is
transverse to $\Sigma$ for all $t$ and the self-intersections of $S_t$
are disjoint from $\Sigma$ for all $t$. We also allow a sphere of $S$
that is disjoint from $\Sigma$ and hence isotopic to a sphere of
$\Sigma$ to be moved from one side of this sphere to the other. The
main result is then:

\begin{thm}\label{proc2.1}
Isotopic systems in normal form are equivalent.\end{thm}

Before beginning the proof of this we will need to establish a more
basic result:

\begin{thm}\label{proc2.2}
Homotopic systems are isotopic. In particular, equivalent systems are
isotopic.\end{thm}

\begin{proof} For single spheres this is Theorem III.1.3 of \cite{L}, and the generalization to systems of spheres follows from the Lemma on page 124 of \cite{L}. We need to extend the proof to allow disks as well. First consider systems consisting of a single disk, so let $ D_0 $ and $ D_1 $ be two disks that are homotopic. If the common boundary of these two disks lies in a torus component $ T $ of $ 
\bdy M $, we can associate to $ D_0 $ and $ D_1 $ spheres $ S_0 $ and
$ S_1 $ that are the boundaries of $\varepsilon${\hy}neighborhoods of
$ T \cup D_0 $ and $ T \cup D_1 $. Alternatively, we can view $ S_i $
as the result of surgering $ T $ along $ D_i $. A homotopy from $ D_0
$ to $ D_1 $ induces a homotopy from $ S_0 $ to $ S_1 $. By
Laudenbach's theorem, $ S_0 $ and $ S_1 $ are then isotopic, so there
is an ambient isotopy of $ M $, fixing $ \bdy M $, taking $ S_0 $ to $
S_1 $ and taking $ D_0 $ to a new disk $ D'_0 $ that lies in the
$\varepsilon${\hy}neighborhood of $ T \cup D_1 $. This neighborhood is
a once-punctured solid torus, the manifold $ M^1_{0,1} $. Any two
disks in $ M^1_{0,1} $ are equivalent under a diffeomorphism of $
M^1_{0,1} $ fixing $ \bdy M^1_{0,1} $ since after filling in the
puncture they are isotopic fixing their boundary. We saw in the
previous section that the mapping class group of $ M^1_{0,1} $ is $
\Z_2 \times \Z $ generated by Dehn twists along the boundary sphere
and boundary torus. Twists along the boundary sphere act trivially on
disks with boundary on the torus. So, up to isotopy, $ D'_0 $ and $
D_1 $ are equivalent under twists along the boundary torus. Looking in
the universal cover of $ M $, lifts of $ D'_0 $ and $ D_1 $ having the
same boundary circle are homotopic and therefore homologous.  But the
homology classes represented by the lifts of the disks obtained by the
$ \Z $'s worth of Dehn twists along $ T $ applied to a single disk are
all distinct, so it follows that $ D'_0 $ and $ D_1 $ are isotopic.

The other possibility is that the common boundary circle of $ D_0 $
and $ D_1 $ lies in a sphere $ S $ of $ \bdy M $. The argument in this
case is similar. An $\varepsilon${\hy}neighborhood of $ S \cup D_i $
is a $3${\hy}punctured sphere bounded by $ S $ and two other spheres $
S'_i $ and $ S''_i $. A homotopy from $ D_0 $ to $ D_1 $ induces a
homotopy from $ S'_0 \cup S''_0 $ to $ S'_1 \cup S''_1 $, so from the
known case of sphere systems we can conclude that these two sphere
systems are isotopic. Arguing as before, this gives a disk $ D'_0 $
isotopic to $ D_0 $ lying in the $\varepsilon${\hy}neighborhood of $ S
\cup D_1 $. Since $ D'_0 $ is not a trivial disk, it is unique up to
isotopy in this neighborhood, which finishes the argument in this
case.

Now we turn to the general case. Let $ X = X_1 \cup \cdotss \cup X_n $
and $ Y = Y_1 \cup \cdotss \cup Y_n $ be two homotopic systems of
spheres and disks, with $ X_i $ corresponding to $ Y_i $ under the
homotopy. By induction on $ n $ we may assume the subsystems $ X_1
\cup \cdotss \cup X_{n-1} $ and $ Y_1 \cup \cdotss \cup Y_{n-1} $ are
isotopic, so after an ambient isotopy we may assume they coincide. The
surfaces $ X_n $ and $ Y_n $ are homotopic in $ M $ so it will suffice
to show they are homotopic in the complement of $ X_1 \cup \cdotss
\cup X_{n-1} $. By another induction the problem reduces further to
the case $ n = 2 $, that is, to show that for systems $ X_1 \cup X_2 $
and $ Y_1 \cup Y_2 $ with $ X_1 = Y_1 $, if $ X_2 $ and $ Y_2 $ are
homotopic in $ M $ then they are homotopic in $ M - X_1 $.

To simplify the notation, let the two systems be $ X \cup Z $ and $ Y
\cup Z $. In the universal cover $ \tild M $ choose a lift $ \tild X $
of $ X $ and lift the homotopy from $ X $ to $ Y $ to a homotopy from
$ \tild X $ to a lift $ \tild Y $ of $ Y $. If $ X $ and $ Y $ are
disks then $ \tild X $ and $ \tild Y $ have the same boundary since
the homotopy from $ X $ to $ Y $ fixes the boundary. After we perturb
$ X $ and $ Y $ to meet transversely we can triangulate $ M $ so that
$ X $, $Y $, and $ Z $ are all subcomplexes, and we can lift this
triangulation to a triangulation of $ \tild M $. Since $ \tild X $ and
$ \tild Y $ are homotopic in $ \tild M $ they are homologous, so there
is a simplicial $3${\hy}chain bounded by $ \tild X $ and $ \tild Y $,
suitably oriented. Geometrically, this chain is a subcomplex $ W $ of
$ \tild M $ bounded by $ \tild X \cup \tild Y $. Note that $ W $ is
unique since there are no $ 3 ${\hy}cycles in $ \tild M $. Also $ W $
is connected since the component of $ W $ containing $ \tild X $ must
also contain $ \tild Y $. This is clear if $ \tild X \cap \tild Y \ne
\emp $, in particular if $ \tild X $ and $ \tild Y $ are disks. If $
\tild X $ and $ \tild Y $ are disjoint spheres they cannot be trivial
in $ H_2(\tild M) \cong \pi_2(\tild M) \cong \pi_2(M) $, so the
component of $ W $ containing $ \tild X $ must contain $ \tild Y $ in
this case as well.

Let $ \tild Z $ be any lift of $ Z $. We claim that $ \tild Z $ is
disjoint from $ W $, except perhaps for points in $ \bdy \tild Z $
when $ \tild Z $ is a disk. For suppose there are other points of
intersection of $ \tild Z $ with $ W $. Then $ \tild Z $ is entirely
contained in $ W $ since $ \bdy W = \tild X \cup \tild Y $ and $ \tild
Z $ is disjoint from $ \tild X $ and $ \tild Y $, apart from $ \bdy
\tild Z $. Since $ \tild Z $ separates $ \tild M $, it splits $ W $
into two parts $ W_1 $ and $ W_2 $. There are only two possibilities:

\items
\item[\sqdot] $ W_1 $ or $ W_2 $ has $ \tild Z $ as its complete
boundary. This is obviously impossible if $ \tild Z $ is a disk, and
if it is a sphere this would make $ \tild Z $ homologous to zero in $
\tild M $, which is not the case since $ Z $ is a nontrivial sphere in
$ M $.

\item[\sqdot] One of $ W_1 $, $ W_2 $ is bounded by $ \tild X \cup
\tild Z $ and the other is bounded by $ \tild Y \cup \tild Z $. In
particular this would say $ \tild X $ and $ \tild Z $ are
homologous. If one is a disk, the other would also have to be a disk
with the same boundary and the sphere $ \tild X \cup \tild Z $ would
be homologous to zero, hence homotopic to zero, making $ \tild X $
homotopic to $ \tild Z $, forcing the same to be true for $ X $ and $
Z $. This would make $ X $ and $ Z $ isotopic, contrary to the
assumption that $ X \cup Z $ is a system. The other alternative is
that both $ \tild X $ and $ \tild Z $ are spheres, but then they would
be homotopic since they are homologous, again forcing $ X $ and $ Z $
to be homotopic and therefore isotopic, a contradiction.  \enditems

Now if we split $ \tild M $ along all the lifts of $ Z $, $ W $ will
lie in one of the resulting components, a simply-connected manifold $
N $. In the case that $ \tild X $ and $ \tild Y $ are spheres, this
says they are homologous in $ N $ and hence homotopic in $ N
$. Projecting to $ M $, we deduce that $ X $ and $ Y $ are homotopic
in the complement of $ Z $, the conclusion we wanted. In the opposite
case that $ \tild X $ and $ \tild Y $ are disks, they together give a
map $ S^2 \rar \tild M $ that is homologous to zero in $ N $, hence
homotopic to zero in $ N $, making $ \tild X $ and $ \tild Y $
homotopic in $ N $. Again this makes $ X $ and $ Y $ homotopic in the
complement of $ Z $. \end{proof}

Before proving Theorem \ref{proc2.1} we need to understand what normal
form systems look like. So suppose the system $ S $ is in normal form
with respect to the maximal sphere system $ \Sigma $. Splitting $ M $
along $ \Sigma $ produces connected manifolds $ P_i $, and splitting $
S $ along $ \Sigma $ produces connected surfaces which we call {\em
pieces} of $ S $. Components of $ S $ disjoint from $ \Sigma $ are
either spheres parallel to components of $ \Sigma $ or disks that are
of one of two types:

\items
\item[(a)] In a $ P_i = S^2 \times I $ that has a component of $
\Sigma $ as one boundary sphere and a component of $ \bdy M $
containing two circles of $ C $ as the other boundary sphere, one can
have a disk with boundary on one circle of $ C $, the disk cutting off
a ball from $ P_i $ that contains the other circle of $ C $.
\item[(b)] In a $ P_i $ that is a punctured $ S^1 \times D^2 $ having
a circle of $ C $ in its boundary torus, one can have a disk with
boundary this circle of $ C $.\enditems

For pieces that actually meet $ \Sigma $ there are the following
possibilities:

\items
\item[(1)] In a $ P_i $ that is a $3 ${\hy}punctured sphere, a piece
can be:

\items
\item[\lowsqdot] a disk with boundary on one boundary sphere and
separating the other two boundary spheres, or
\item[\lowsqdot] a cylinder with its two boundary circles on two
different boundary spheres of $ P_i $, or
\item[\lowsqdot] a pair of pants with each boundary circle on a
different boundary sphere of $ P_i $.  \enditems

\item[(2)] In a $ P_i = S^2 \times I $ with a component of $ \Sigma $
at one end and a component of $ \bdy M $ at the other, a piece can be
a cylinder with one boundary circle on each boundary sphere, the
boundary circle in $ \bdy M $ being a circle of $ C $.

\item[(3)] In a $ P_i $ that is a once-punctured summand $ Q_j $, a
piece can be the boundary surface of a tubular neighborhood of a tree
in $ P_i $ obtained by joining $ p \ge 2 $ points in the boundary
sphere of $ P_i $ to an interior point of $ P_i $ by disjoint arcs
such that a lift of this tree to the universal cover $ \tild P_i $ of
$ P_i $ has its endpoints all in distinct boundary spheres of $ \tild
P_i $. This is equivalent to saying that the boundary circles of a
lift of the piece all lie in different components of $ \bdy \tild P_i
$.

\item[(4)] In a $ P_i $ that is a once-punctured $ S^1 \times D^2 $
with a circle of $ C $ in its boundary torus a piece can be of the
type in (3) or it can be:

\items
\item[\lowsqdot] a cylinder joining the circle of $ C $ in the
boundary torus to the boundary sphere, or
\item[\lowsqdot] a pair of pants with one boundary circle on $ C $ and
two boundary circles on the boundary sphere of $ P_i $, such that
after lifting to the universal cover the three boundary circles of a
lift of the pair of pants all lie in different components of $ \bdy
\tild P_i $.  \enditems

\enditems

\begin{figure}[ht!]\small\anchor{fig3}
\SetLabels (.25*.21) $\Sigma$\\ (.407*.85) $S_0$\\ \endSetLabels
\centerline{\AffixLabels{\epsfbox{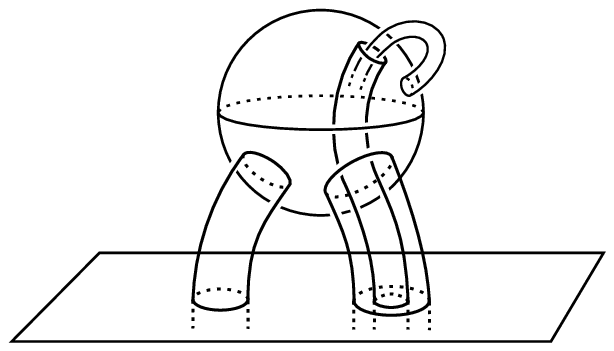}}}
\nocolon\caption{}\label{fig3}
\end{figure}

To show that if the number of circles of intersection of $ S $ with $
\Sigma $ is minimal then the pieces have only the types listed we
argue as follows. By a sequence of surgeries on $ S $ the circles of $
S \cap \Sigma $ can be eliminated one by one. This converts $ S $ into
a collections of disks and spheres disjoint from $ \Sigma $, and $ S $
is obtained from this collection by the inverse sequence of tubing
operations, where the tubes might be nested one inside another. So
each piece of $ S $ is obtained from a sphere or disk in a $ P_i $ by
inserting tubes that run from this sphere or disk to one or more of
the boundary spheres of $ P_i $ that are in $ \Sigma $. See
\figref{fig3} for an example.

Consider first the case of a sphere $ S_0 $ connected to boundary
spheres of $ P_i $ by tubes. Since $ \Sigma $ is maximal, $ S_0 $ must
be trivial in $ P_i $, either bounding a ball or parallel to a sphere
of $ \bdy P_i $. Suppose first that $ S_0 $ bounds a ball. If there
were nesting among the tubes from $ S_0 $, then there would be an
outermost tube connecting $ S_0 $ on the outside to a sphere of $
\Sigma $ and a next-outermost tube inside this going from the inside
of $ S_0 $ to the same sphere of $ \Sigma $.

There is an isotopy of $ S $ that takes the end of this inner tube
attached to $ S_0 $ and slides it along $S_0$ over to the outer tube
and then along the outer tube until it lies on the other side of $
\Sigma $. If there are other tubes inside this inner tube, they are to
be carried along during this sliding process. After this isotopy the
number of circles of $ S \cap \Sigma $ has temporarily increased, but
now there is a regular homotopy of the new tube that pushes it across
$ \Sigma $ so as to eliminate both its new and old intersections with
$ \Sigma $. Again any tubes inside this tube can be carried along
during the regular homotopy. The net result is a new embedding of $S$
which is homotopic to the original one and therefore isotopic to it by
Theorem~\ref{proc2.2}, and which has fewer circles of intersection
with $\Sigma$ than the original $S$ had, contrary to the minimality
assumption. Thus nesting of tubes is ruled out. Similarly, if two
tubes joining $ S_0 $ to the same sphere of $ \Sigma $ were homotopic
in $ P_i $, then we could slide one tube across the other and again
decrease the number of circles of $ S \cap \Sigma $. Thus homotopic
tubes are ruled out, or in terms of the universal cover, tubes joining
a lift of $ S_0 $ to the same lift of a sphere of $ \Sigma $. Since
there have to be at least two tubes from $ S_0 $, we see in this case
that the piece is as in (3) or is a cylinder or pair of pants as in
(1).

Suppose now that $ S_0 $ is a sphere parallel to a sphere of $ \bdy
P_i $. If there is a tube from $ S_0 $ to this sphere of $ \bdy P_i $,
then the piece can also be obtained by joining a sphere bounding a
ball to $ \bdy P_i $ by tubes, reducing to the previous case. If there
is no tube joining $ S_0 $ to the sphere of $ \bdy P_i $ it is
parallel to, then $ P_i $ must be a $3${\hy}punctured sphere as in
(1). By the reasoning in the preceding paragraph there can be no
nesting of tubes, and there can be at most one tube joining $ S_0 $ to
each of the other two spheres of $ \bdy P_i $ since $ P_i $ is
simply-connected. Thus the piece is a disk or a cylinder, and it must
be of the type in (1).

The other possibility is that the piece is a disk $ D_0 $ with tubes
connecting it to $ \Sigma $. If $ P_i $ is as in (2), then arguing as
before we see that there can be only one tube and the piece is a
cylinder of the type described. If $ P_i $ is as in (4) then $ D_0 $
is a nontrivial disk in $ P_i $ joined by one or more tubes to the
boundary sphere of $ P_i $. The disk $ D_0 $ is unique up to isotopy
if we allow its boundary to move in the boundary torus of $ P_i $, so
we may assume $ D_0 $ is obtained by applying Dehn twists along this
torus to a meridian disk of $ P_i $. Minimality rules out nesting of
tubes, and it rules out two tubes being connected to the same side of
$ D_0 $ since if there were, the two tubes would be homotopic. Thus
either there is one tube or there are two tubes attached on opposite
sides of $ D_0 $. These are the configurations in (4).

\begin{proof}[Proof of Theorem \ref{proc2.1}]  What we will show is that if two isotopic systems $ S $ and $ S'
$ are in normal form in the sense that their pieces have the forms
described above, then they are equivalent.

Let $ \tild \Sigma $ be the preimage of $ \Sigma $ in the universal
cover $ \tild M $ of $ M $. We first consider the case that $ S $ has
a single component, a sphere or disk. Let $ \tild S $ be a lift of $ S
$ to $ \tild M $. An isotopy of $ S $ to $ S' $ lifts to an isotopy of
$ \tild S $ to a lift $ \tild S' $ of $ S' $. Our goal in the first
part of the proof will be to show:

\items
\item[{\bf (I)}] $ \tild S $ and $ \tild S' $ are isotopic staying
transverse to $ \tild\Sigma $, except in the special case that they
are spheres isotopic to the same sphere of $ \tild\Sigma$ but lying on
opposite sides of this sphere.  \enditems

The proof will involve several steps. To begin, let us choose an
orientation for $ \tild S $. This is equivalent to choosing a
transverse orientation once we fix an orientation for $ \tild M $. The
orientation for $ \tild S $ carries over to an orientation for $ \tild
S' $ by means of the isotopy between them, and likewise for the
transverse orientation.

Dual to $ \tild \Sigma $ is a tree $ T $, with a vertex for each
component of $ \tild M - \tild \Sigma $ and an edge for each component
of $ \tild \Sigma $.  Similarly, dual to the collection of circles $
\tild S \cap \tild\Sigma $ in $ \tild S $ is a tree $ T(\tild S) $
with a vertex for each component of $ \tild S - \tild \Sigma $ and an
edge for each circle of $ \tild S \cap \tild\Sigma $. The natural map
$ T(\tild S) \rar T $ is injective since it is locally injective, the
pieces of $ S $ having the types listed in (1)--(4) when $ S \cap
\Sigma \ne \emp $. Thus we can view $ T(\tild S) $ as a subtree of $ T
$. In the same way we have a subtree $ T(\tild S') $ of $ T $.

\items
\item[{\bf (a)}] $ T(\tild S) = T(\tild S') $, except in the case that
$ S $ is a sphere isotopic to a sphere of $ \Sigma $, when $ T(\tild
S) $ and $ T(\tild S') $ can be adjacent vertices of $ T $.  \enditems

To prove this it will suffice to show that $ T(\tild S) $ is
determined by the homology class of the oriented surface $ \tild S $
in $ H_2(\tild M,\bdy\tild S) $.  Let $e$ be an oriented edge of $ T
$, by which we mean an edge together with a choice of one of its two
possible orientations. We can split $ \tild M $ into two components
along the sphere of $ \tild\Sigma $ corresponding to $ e $, and we let
$ \tild M_e $ be the component at the tail of~$ e $. If $ T(\tild S) $
is not a single vertex, then an end vertex $ v $ of $ T(\tild S) $
corresponds to a piece of $ \tild S $ that is either a disk in a
$3${\hy}punctured sphere or an annulus in a once-punctured $D^3 $ or $
S^1 \times D^2 $.  This disk or annulus is nonzero in $ H_2(\tild M,
\tild M_e \cup \bdy\tild S) $ where $ e $ is the edge of $ T(\tild S)
$ with $ v $ as its head vertex. Thus if we partially order the
submanifolds $ \tild M_e $ by inclusion as $ e $ ranges over the
oriented edges of $ T $, then the maximal $ \tild M_e $'s for which $
\tild S $ is nonzero in $ H_2(\tild M,\tild M_e \cup \bdy\tild S) $
are those corresponding to oriented edges $ e $ whose head vertex is
an end vertex of $ T(\tild S) $, at least if $ T(\tild S) $ is not a
single vertex. This remains true also when $ T(\tild S) $ is a single
vertex and $ S $ is a disk. In the case that $ T(\tild S) $ is a
single vertex and $ S $ is a sphere the maximal $ \tild M_e $'s have $
e $ with its head vertex at either end of the edge of $ T $
corresponding to the component of $ \tild\Sigma $ isotopic to $ \tild
S $. This gives a homological characterization of $ T(\tild S) $,
proving (a).

Assuming that we are not in the exceptional case in (a) from now on,
then since $ T(\tild S) = T(\tild S') $, the pieces of $ \tild S $ are
diffeomorphic to those of $ \tild S' $. In fact a stronger statement
is true:

\items
\item[{\bf (b)}] The pieces of $ \tild S $ are isotopic within their $
\tild P_i $'s to the corresponding pieces of $ \tild S' $.  Here we
mean proper isotopies of the pieces, so their boundary circles in $
\tild\Sigma $ stay in $ \tild\Sigma $ during the isotopies.  \enditems

This is proved by examining the various types of pieces. If a piece is
the boundary of a neighborhood of a tree in $ \tild P_i $ then two
different embeddings of the tree having the same endpoints are
homotopic fixing their endpoints since $ \tild P_i $ is
simply-connected, and a homotopy can be improved to an isotopy by
sliding crossings of edges over the boundary spheres at their outer
ends --- the so-called lightbulb trick of unknotting the cord to a
hanging light bulb by letting the knot drop over the bulb. This works
in the present situation since the various edges in the tree all go to
different boundary spheres of $ \tild P_i $. This argument covers
pieces of types (2) and (3) and cylinders and pairs of pants in (1). A
disk in (1) is unique up to isotopy since the isotopy class of a disk
in a ball with punctures is determined by how it separates the
punctures. For (4) we have a piece that is the lift of a meridian disk
with one or two tubes attached, and possibly with Dehn twists along
the boundary torus applied to the meridian disk. Again the lightbulb
trick applies in $ \tild P_i $ since if there are two tubes, they
connect to different boundary spheres of $ \tild P_i $. Thus we have
(b) in all cases.

The different ways of gluing together two pieces that lie on opposite
sides of a sphere of $ \tild\Sigma $ and have a boundary circle in
this sphere depend only on choosing an isotopy taking one boundary
circle to the other in this sphere. The space of circles in a $ 2
${\hy}sphere has the homotopy type of the projective plane, with
fundamental group $ \Z_2 $, so there are two essentially distinct
isotopies taking one circle to another. This two-fold ambiguity is
eliminated if we require the gluing to preserve transverse
orientations for the two pieces since the space of oriented circles in
a $2${\hy}sphere is simply-connected, being homotopy equivalent to $
S^2 $. Thus in order to finish the proof of assertion (I) that $ \tild
S $ and $\tild S' $ are isotopic staying transverse to $ \tild\Sigma
$, it will suffice to refine (b) to:
 
\items
\item[{\bf (c)}] Corresponding pieces of $ \tild S $ and $ \tild S' $
are isotopic preserving transverse orientations.  \enditems

This is automatically true for a piece that is a pair of pants in a
$3${\hy}punctured sphere since there is an isotopy of this piece that
starts and ends with the same position and reverses the transverse
orientation. Namely, take the three punctures to lie on a great circle
in $ S^3 $ and then the family of great $ 2 ${\hy}spheres containing
this great circle gives such an isotopy. (Alternatively, think of the
$3${\hy}punctured sphere as a $2${\hy}punctured ball, with the two
punctures aligned along a diameter of the ball, and intersect the
punctured ball with the family of planes containing this diameter.)
The same thing happens for a piece that is an $ n ${\hy}punctured
sphere in a $ \tild P_i $ that is an $ n ${\hy}punctured $ S^3 $. Such
a $ \tild P_i $ is the universal cover of a once-punctured $ Q_j $
with $ \pi_1 Q_j $ of order $ n $.

We claim that cutting $ \tild S $ along $ \tild\Sigma $ and regluing
according to different transverse orientations in any subset of the
remaining pieces changes the class of $ \tild S $ in $ H_2(\tild
M,\bdy\tild S) $. If $ \bdy S \ne \emp $ it is clear that reorienting
the piece containing $ \bdy S $ changes the homology class since this
changes its image under the boundary map to $ H_1(\bdy\tild S) $.  So
we may assume the orientation of this piece is fixed. To treat the
remaining cases it is convenient to pass to a quotient manifold $ N $
of $ \tild M $ constructed in the following way. Consider a $ \tild
P_i $ containing a piece of type (3) that does not touch all the
boundary spheres of $ \tild P_i $. This piece is the boundary of a
neighborhood $ X $ of a tree in $ \tild P_i $. Let $ Y $ be the union
of $ X $ with an arc joining it to a boundary sphere of $ \tild P_i $
disjoint from the piece. Let $ Z $ be an $ \epsilon ${\hy}neighborhood
in $ \tild P_i $ of the union of $ Y $ with the boundary spheres of $
\tild P_i $ that meet $ Y $. Then $ \bdy Z $ consists of some boundary
spheres of $ \tild P_i $ together with one sphere in the interior of $
\tild P_i $. The latter sphere splits $ \tild M $ into two
components. Collapsing the component not containing $ Z $ to a point
produces a quotient manifold of $ \tild M $. Doing this collapsing
operation on all such $ \tild P_i $'s gives the manifold $ N $. This
still contains $ \tild S $, and it contains a subcollection $ \tild
\Sigma_N $ of the spheres of $ \tild\Sigma $, with dual tree $ T_N $
containing $ T(\tild S) $. Edges of $ T_N $ that touch $ T(\tild S) $
but are not contained in $ T(\tild S) $ we call {\em abutting} edges.

Consider the effect of changing the transverse orientation on a type
(3) piece that meets all but one of the boundary spheres of the
$n${\hy}punctured sphere of $ N - \tild\Sigma_N $ that contains it. If
we fill in with a ball the one boundary sphere that it does not meet,
then we have seen that there is an isotopy of the piece that reverses
its transverse orientation.  This isotopy sweeps across the filled-in
ball exactly once, so if we take out the filled-in ball we see that
reorienting the piece changes the homology class of $ \tild S $ by
adding the boundary sphere of the filled-in ball, with one of its two
possible orientations. The same argument applies also to a type (1)
piece that is an annulus. For a type (1) piece that is a disk, if we
surger $ \tild S $ along the boundary of this disk using one of the
two disks it bounds in $ \tild\Sigma $, we obtain a surface of two
components, one of which is a copy of one of the other two boundary
spheres of the $ \tild P_i $ containing the given piece, and we can
recover $ S $ from this surgered surface by a connected sum operation
reversing the surgery. The connected sum operation preserves the
homology class. As \figref{fig4} shows, changing the boundary sphere
we use for the connected sum corresponds to changing the orientation
of the piece.

\begin{figure}[ht!]\anchor{fig4}
\centerline{\epsfbox{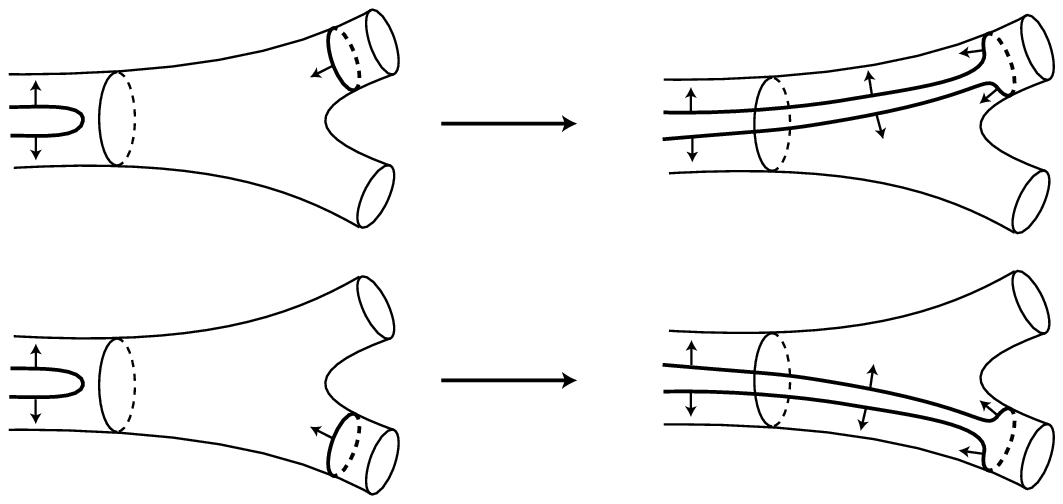}} \nocolon\caption{}\label{fig4}
\end{figure}

Hence the homology class of $ \tild S $ changes by subtracting one
boundary sphere of $ \tild P_i $ and adding another. In terms of
transverse orientations the homology class changes by adding two
boundary spheres of $ \tild P_i $, one oriented into $\tild P_i $ and
the other out of $ \tild P_i $.

Thus we see that reorienting some collection of pieces disjoint from $
\bdy\tild S $ changes the homology class of $ \tild S $ in $
H_2(N,\bdy\tild S) $ by adding the spheres corresponding to some of
the abutting edges of $ T_N $, with certain signs. If this signed sum
of spheres were homologically a boundary it would have to consist of
spheres for all the abutting edges, oriented all toward $ \tild S $ or
all away from $ \tild S $, but this is not the case as there is at
least one abutting edge corresponding to a disk. This finishes the
proof of (c) and hence of (I).

We now consider the general case that $ S $ and $ S' $ have an
arbitrary number of components. Choose a transverse orientation for $
S $. This gives a transverse orientation to the circles of $ S \cap
\Sigma $ in $ \Sigma $. In each sphere $ \Sigma_0 $ of $ \Sigma $,
each circle of $ S \cap \Sigma_0 $ then partitions the remaining
circles of $ S \cap \Sigma_0 $ into those on the positive side and
those on the negative side of the given circle. Let us call this
partition data the {\em side relation\/} for $ S \cap \Sigma $. The
transverse orientation for $ S $ induces a transverse orientation for
$ S' $ via the isotopy between them, so the circles of $ S' \cap
\Sigma $ also have a side relation. From what we have shown so far we
have a bijective correspondence between the pieces of $ S $ and $ S' $
and hence also between the circles of $ S \cap \Sigma $ and those of $
S' \cap \Sigma $.

\items
\item[{\bf (II)}] The side relation on $ S \cap \Sigma $ agrees with
that on $ S' \cap \Sigma $ under this bijection.  \enditems

For suppose this is false. Then in some sphere $ \Sigma_0 $ of $
\Sigma $ we have circles $ c_p $ and $ c_q $ of $ S \cap \Sigma_0 $
with corresponding circles $ c'_p $ and $ c'_q $ of $ S' \cap \Sigma_0
$ such that $ c_q $ and $ c'_q $ lie on opposite sides of $ c_p $ and
$ c'_p $, respectively. Choose a lift $ \tild \Sigma_0 $ of $ \Sigma_0
$ to $ \tild M $, and let $ \tild S_j $ and $ \tild S_k $ be the lifts
of components of $ S $ containing the lifts to $ \tild \Sigma_0 $ of $
c_p $ and $ c_q $, respectively. Similarly we have $ \tild S'_j $ and
$ \tild S'_k $ containing the lifts of $ c'_p $ and $ c'_q $ to $
\tild\Sigma_0 $. We know that $ \tild S'_j $ is isotopic to $ \tild
S_j $ staying transverse to $ \tild\Sigma $. This isotopy can be
realized by an ambient isotopy of $ (\tild
M,\tild\Sigma)$. Restricting this ambient isotopy to $ \tild S'_k $,
we obtain an isotopy of $ \tild S'_k $ to a surface $ \tild S^*_k
$. Then $ \tild S_k $ and $\tild S^*_k $ lie on opposite sides of
$\tild S_j $, hence are disjoint except for their common boundary
circle when they are disks. Since they are isotopic, both being
isotopic to $ \tild S'_k $, there is a product region between
them. This product contains $ \tild S_j $ since $ \tild S_k $ and $
\tild S^*_k $ lie on opposite sides of $ \tild S_j $. Thus if $\tild
S_k $ and $ \tild S^*_k$ are spheres, so is $S_j$, and if they are
disks, so is $S_j$. There is a unique isotopy class of disks or
nontrivial spheres in the product between $ \tild S_k$ and $\tild
S^*_k$, so there is also a product region between $ \tild S_j $ and $
\tild S_k $. If $ \tild S_j $ and $ \tild S_k $ project to different
components of $ S $, these two components would be homotopic and hence
isotopic, contrary to the definition of a system of disks and
spheres. So $ \tild S_j $ and $ \tild S_k $ must be lifts of the same
component of $ S $. This component must be a sphere, otherwise the
isotopy between $ \tild S_j $ and $ \tild S_k $ would give $ \bdy
\tild S_j = \bdy \tild S_k $ and hence $ \tild S_j = \tild S_k $,
which is impossible since $ c_p \ne c_q $.

\begin{lem}\label{proc2.3}
If a $3${\hy}manifold $M$ contains a nontrivial sphere $S$ having two
lifts to the universal cover $\tild M$ that bound a product $S^2
\times I$, then $ M $ is either an $ S^2 ${\hy}bundle over $S^1$ or
$M$ has a connected summand that is a real projective $3${\hy}space
split off by $S$.\end{lem}

\begin{proof} Let $S_1$ and $S_2$ be two lifts of $S$ bounding a product $S^2 \times I$ in $
\tild M$. By rechoosing $S_2$ if necessary, we may assume no other
lifts of $S$ lie in this product.  If the deck transformation taking
$S_1$ to $S_2$ takes the product to itself, it is a free involution on
this product, yielding a projective space summand of $ M $ bounded by
$S$ by a classical theorem of Livesay \cite{Li}. If the deck
transformation takes the product $ S^2 \times I $ outside itself, then
$M$ is obtained from this product by identifying its two ends, so $M$
is an $ S^2${\hy}bundle. \end{proof}

Now we finish the proof of (II). If $M$ is an $ S^2${\hy}bundle, it
contains a unique isotopy class of nontrivial spheres, and if $S$ and
$S'$ were in normal form with respect to $\Sigma$ they would have to
be disjoint from $\Sigma$, making (II) true vacuously.  So we may
assume $M$ is not an $ S^2${\hy}bundle. In the case of a projective
space summand bounded by $S_j$, this summand lifts to the region
between $ \tild S_j $ and $ \tild S_k $. Applying the same argument
with the roles of $ S $ and $ S' $ reversed, we find a projective
space summand of $M$ bounded by $S'_j$ lifting to the region between $
\tild S'_j$ and $ \tild S'_k$. After extending the isotopy taking
$S'_j$ to $S_j$ to an isotopy of $M$, the projective space summand
bounded by $S'_j$ becomes a second projective space summand bounded by
$S_j$, and this lies on the opposite side of $S_j$ from the first
summand since $\tild S_k $ and $ \tild S'_k $ lie on opposite sides of
$\tild S_j $ and $ \tild S'_j $, where we are using the transverse
orientations to distinguish sides. Thus $S_j$ has a projective space
summand on each side, so $M$ is the connected sum of two projective
spaces.  It is an elementary fact that in a manifold such as this
which is the sum of two irreducible manifolds there is only one
isotopy class of nontrivial spheres, just as for $ S^2${\hy}bundles,
and (II) would again be true vacuously. (If one wants to avoid quoting
Livesay's theorem, one can say that $S_j$ splits $M$ as the sum of two
closed $3${\hy}manifolds with fundamental group of order $2$, and
these manifolds are irreducible since $M$ has no fake $3${\hy}sphere
summands.)

\rk {Conclusion of the proof of Theorem \ref{proc2.1}} It is easy to
see that a system of transversely oriented circles in a sphere is
determined up to isotopy by its associated side relation. In view of
(II) this means that by an isotopy of $ S $ staying transverse to $
\Sigma $, the transversely oriented circles of $ S \cap \Sigma $ can
be made to agree with the corresponding transversely oriented circles
of $ S' \cap\Sigma $. By the earlier arguments there are then
isotopies of the pieces of $ \tild S $ to the corresponding pieces of
$ \tild S' $, preserving transverse orientations, and we can take
these isotopies to fix the boundaries of these pieces since the space
of transversely oriented circles in a $2${\hy}sphere is
simply-connected, as noted earlier. Thus these isotopies fit together
to give an isotopy of $ \tild S $ to $ \tild S' $ fixing $ \tild S
\cap \tild \Sigma $. This isotopy projects to an equivalence of $ S $
and $ S' $ in $ M $. \end{proof}

\section{Complexes of disks and spheres}

We continue with the notation of the previous section, so $ M $ is a
connected sum of $ n $ copies of $ S^1 \times S^2 $ with $ s $ copies
of the ball $ D^3 $ and $ k $ additional irreducible manifolds, and $
C $ is a collection of circles in the sphere and torus boundary
components of $ M $ subject to the same restrictions as before. Let $
\DS(M,C) $ be the simplicial complex whose $p${\hy}simplices are the
isotopy classes of systems of ${p+1}$ disks and spheres in $ M $. The
faces of such a simplex are obtained by passing to subsystems. If we
restrict attention to systems formed just of disks or just of spheres
we obtain subcomplexes $ \D(M,C) $ and $ \Sph(M,C)$, respectively.  We
write $ \Sph(M) $ for $ \Sph(M,\emp) $.  Eventually we will be most
interested in the subcomplex $\D_c(M,C)$ of $\D(M,C)$ formed by disk
systems with connected complement. The main result we are heading for
in this section is

\begin{thm}\label{proc3.1}
$\D_c(M,C)$ is $(n-2)${\hy}connected if $ C \ne \emp $, and
$(n-1)${\hy}connected if $C$ is contained in torus components of $\bdy
M$ and $ C \ne \emp $. \end{thm}

The plan of the proof is to show:

\items
\item[(1)] $\DS(M,C)$ is contractible except when $ n=0 $ and $ k \le
1 $, in which cases it is ${(2k+s-5)}${\hy}connected if $ C = \emp $
and $(k+s-4)${\hy}connected if $C \ne \emp$.

\item[(2)] $\D(M,C)$ is contractible if $ C $ is contained in torus
components of $ \bdy M $ and $ C \ne \emp$.

\item[(3)] $\D(M,C) $ is $(2n+k+s-4)${\hy}connected in general, if $ C
\ne \emp$.

\item[(4)] $\D_c(M,C)$ is $(n-2)${\hy}connected if $ C \ne \emp$.

\item[(5)] $\D_c(M,C)$ is $(n-1)${\hy}connected if $C$ is contained in
torus components of $\bdy M$, $ C \ne \emp$.  \enditems

Our convention is that $(-1)${\hy}connected means nonempty and
$(-2)${\hy}connected is a vacuous condition. The empty set is
$(-2)${\hy}connected, for example.

The proof of (3) is the most delicate so we will postpone this until
last.

\begin{proof}[Proof of (1)] There are three main steps. The first is the following, where $ t $ denotes the number of boundary spheres of $ M $ that contain circles of $ C $:

\items
\item[{\bf (a)}] $\DS(M,C) $ is contractible if $ \Sph(M) \ne \emp $
and either $ s \le 1 $ or $ s = t $.  \enditems

This is proved using the surgery technique of the proof of Theorem~2.1
of \cite{H3}, with a few minor simplifications. Choose a maximal
sphere system $\Sigma$ and put an arbitrary disk-sphere system $ S =
S_0 \cup \cdotss \cup S_p $ in normal form with respect to~$ \Sigma
$. A point in the $p${\hy}simplex of $ \DS(M,C) $ determined by $ S $
can be thought of as a weighted sum $ \sum_i t_iS_i $ of the
components of $ S $, with weights $t_i $ the barycentric coordinates
of the point in the simplex. Barycentric coordinates are normalized to
have $ \sum_i t_i = 1 $. To make the surgery process clearer it is
helpful to replace each $ S_i $ by a family $ S \times [0,t_i] $ of
parallel copies of $ S_i $ of `thickness' $ t_i $. When a weight $ t_i
$ goes to zero at a face of the simplex, the family $ S_i \times
[0,t_i] $ shrinks to thickness zero and is deleted. We will allow the
family $ S_i \times [0,t_i] $ to be split into several parallel
families of total thickness $ t_i $, as well as the inverse operation
of combining parallel families into a single family, adding their
weights.

Our aim is to construct a sequence of surgeries on $ S $ to eliminate
the circles of $ S \cap \Sigma_0 $ for $ \Sigma_0 $ one of the spheres
of $ \Sigma $. Let $ T_0 $ be the dual tree of $ S \cap \Sigma_0 $ in
$ \Sigma_0 $, with a vertex for each components of $ \Sigma_0 - S $
and an edge for each circle of $ S \cap \Sigma_0 $. The weights on the
$ S_i $'s define lengths for these edges, so $ T_0 $ is a metric
tree. There is a canonical way to shrink $ T_0 $ to a point by
shortening all its extremal edges simultaneously at unit speed. Once
an extremal edge has disappeared one continues shrinking all remaining
extremal edges. The surgery process will realize this shrinking of $
T_0 $ by a path in $ \DS(M,C) $ starting at $ S $ and ending with a
system disjoint from $ \Sigma_0 $. If $ T_0 $ is not already a point,
its extremal vertices $ v $ correspond to disjoint disks $ D_v $ in $
\Sigma_0 $ with $ D_v \cap S = \bdy D_v $. These disks can be used to
surger $ S $ to a new system in which the circles $ \bdy D_v $ have
been eliminated from $ S \cap \Sigma_0 $.  Taking weights into
account, one gradually surgers through the appropriate families $ S_i
\times [0,t_i] $ at unit speed, decreasing the thicknesses of these
families while increasing the thicknesses of the new families created
by the surgery. The old and new families can be taken to be disjoint,
so the surgery can be viewed as simply transferring weights from the
old family to the new. When the thickness of a family has shrunk to
zero, this family is deleted and one continues the surgery process on
the remaining families.  Near the end of the process of shrinking $
T_0 $ to a point it can happen that all that remains of $ T_0 $ is a
single edge, and then both ends of this edge are shrinking at the same
time, so the corresponding family of surfaces is being surgered on
both sides simultaneously.

Thus the surgery process produces a family $ S(t_1,\cdotss,t_p,u) $ of
weighted collections of surfaces, for $ (t,u) \in \Delta^p \times
[0,\infty) $, with $ S(t_1,\cdotss,t_p,0) $ the given family. After
the process is completed we combine parallel families and discard any
trivial spheres or disks that are produced by the surgeries. We need
to make sure that some nontrivial disks or spheres always remain.
Surgery on a nontrivial sphere $ S_0 $ produces a pair of spheres, and
the only way these can both be trivial is if $ S_0 $ splits off a pair
of boundary spheres from $ M $, neither of which meets $ C $.  This
cannot happen if $ s \le 1 $ or $ s = t $. Surgery on a nontrivial
disk $ D_0 $ produces a sphere and a disk, and the only way these can
both be trivial is if $ \bdy D_0 $ lies in a sphere of $ \bdy M $ and
$ D_0 $ splits off another boundary sphere of $ M $ that is disjoint
from $ C $. Again we are safe if $ s \le 1 $ or $ s = t $.

Thus after discarding trivial spheres and disks we are left with
systems of spheres and disks.  Then renormalizing the weights of the
remaining nontrivial spheres and disks to have sum $ 1 $, the surgery
process defines a path in $ \DS(M,C) $. This is independent of the
choice of normal form system $ S $ within its isotopy class because
the surgeries always take equivalent systems to equivalent systems,
and equivalent systems are isotopic. The surgeries may destroy normal
form, but that does not matter.

It is evident that $ S(t_1,\cdotss,t_p,u) $, regarded as a point in $
\DS(M,C) $ by discarding trivial components and renormalizing weights,
depends continuously on the parameters $ t_i $ and $ u $. When we pass
to a face of the parameter simplex $ \Delta^p $ by letting some $ t_i
$ go to zero, we get the corresponding family $ S(t_1,\cdotss,\widehat
t_i,\cdotss,t_p,u) $ for the system obtained by deleting $ S_i $ from
$ S $. So we have constructed a deformation retraction of $ \DS(M,C) $
into the subspace of systems disjoint from $ \Sigma_0 $. This subspace
is contractible since it is the star of the vertex corresponding to $
\Sigma_0 $. Hence $ \DS(M,C) $ is contractible.

Note that instead of surgering along just the sphere $ \Sigma_0 $ we
could just as well surger along any subsystem of $ \Sigma $, and this
would produce a deformation retraction onto the star of this subsystem
in $ \DS(M,C) $.

This takes care of (a). The other two steps in the proof of (1) will
provide a reduction to (a). Suppose $ M $ has boundary spheres $
\bdy_1,\cdotss,\bdy_s$, one of which, say $ \bdy_s $, is disjoint from
$ C $. Let $ M' $ be the manifold obtained from $ M $ by filling in $
\bdy_s $ with a ball.

\items
\item[{\bf (b)}] If $ \Sph(M) \ne \emp $ and $ s \ge t+2 $, the
connectivity of $ \DS(M,C) $ is one greater than the connectivity of $
\DS(M',C) $.  \enditems

The proof of this is similar to the proof of Lemma~2.2 in \cite{H3}.
Call a vertex of $ \DS(M,C) $ {\em special\/} if it is either

\items
\item[\lowsqdot] a sphere splitting off a $3${\hy}punctured sphere
from $ M $ whose other two boundary spheres are $ \bdy_s $ and some $
\bdy_i $ which is disjoint from $ C $, or
\item[\lowsqdot] a disk splitting off a $2${\hy}punctured sphere from
$ M $ having $ \bdy_s $ as one boundary sphere and whose other
boundary sphere intersects $ C $ only in the boundary circle of the
disk.  \enditems

Special vertex spheres exist if $ s \ge t+2 $, so let $ \Sigma $ be
such a sphere, splitting off a submanifold $ P $ of $ M $ that is a $
3 ${\hy}punctured sphere bounded by $ \Sigma $, $ \bdy_s $, and some
other $ \bdy_i $. Let $ \DS'(M,C) $ be the subcomplex of $ \DS(M,C) $
consisting of simplices with no special vertices. A deformation
retraction of $ \DS'(M,C) $ onto the link of $ \Sigma $ can be
obtained in the following way. Take a system $ S $ defining a simplex
in $ \DS'(M,C) $ and put it in normal form with respect to some
maximal sphere system containing $ \Sigma $.  Then if $ S \cap P $ is
not empty, it consists of parallel disks separating $ \bdy_s $ from $
\bdy_i $. These disks can be eliminated by pushing them across $
\bdy_s $ and then outside $ P $. Since no components of $ S $ are
special disks or spheres, this process of modifying $ S $ by pushing
across $ \bdy_s $ produces no trivial disks or spheres. The process
could also be described in terms of surgering $ S $ along $ \Sigma $
and discarding the resulting spheres that are isotopic to $ \bdy_s $,
so the process determines a deformation retraction of $ \DS'(M,C) $
onto the link of $ \Sigma $, as desired.

We obtain $ \DS(M,C) $ from $ \DS'(M,C) $ by attaching the star of
each special vertex along the link of this vertex. These stars are the
cones on the links, and have disjoint interiors since no simplex can
contain two distinct special vertices.  The star of $ \Sigma $ is
contractible, so by the preceding paragraph the union of this star
with $ \DS'(M,C) $ is also contractible. We obtain $ \DS(M,C) $ from
this contractible space by attaching the stars of the other special
vertices, so $ \DS(M,C) $ is homotopy equivalent to the wedge of the
suspensions of the other links. Each link is a copy of $ \DS(M',C) $,
so the connectivity of $ \DS(M,C) $ is one greater than the
connectivity of $ \DS(M',C) $.

\items
\item[{\bf (c)}] The connectivity of $ \DS(M,C) $ is one greater than
the connectivity of $ \DS(M',C) $ if $ s > t $ and $ t > 0 $. The same
is true for $ \D(M,C) $ and $ \D(M',C) $.  \enditems

The proof parallels the previous one. First we have a preliminary
observation:

\items
\item[\sqdot] If $ C $ is obtained from $ C' $ by adding a second
circle in some sphere component of $ \bdy M $ then $ \DS(M,C) $ is
homotopy equivalent to the suspension of $ \DS(M,C') $, and this is
also true for $ \D(M,C) $ and $ \D(M,C') $.  \enditems

To see this, let $ C_0 $ and $ C_1 $ be two circles of $ C $ in the
same boundary sphere of $ M $ and let these bound disks $ D_0 $ and $
D_1 $ that lie in a neighborhood of the boundary sphere and represent
vertices of $ \D(M,C) $. Let $ \DS'(M,C) $ be the subcomplex of $
\DS(M,C) $ consisting of simplices containing neither of the vertices
represented by $ D_0 $ and $ D_1 $. Then $ \DS'(M,C) $ deformation
retracts onto the link of $ D_1 $ in $ \DS(M,C) $ by shifting disks
with boundary on $ C_0 $ over to disks with boundary on $ C_1 $, using
the annulus between $ C_0 $ and $ C_1 $ in $ \bdy M $.  In the same
way $ \DS'(M,C) $ deformation retracts onto the link of $ D_0 $. Since
$ \DS(M,C) $ is the union of $ \DS'(M,C) $ with the stars of the
vertices $ D_0 $ and $ D_1 $, it follows that $ \DS(M,C) $ is homotopy
equivalent to the suspension of the link of either of these
vertices. These links can be identified with $ \DS(M,C') $. The same
argument applies with $ \D $ in place of $ \DS $.

Now to prove (c) we may assume $ C $ has at most one circle in each
sphere of $ \bdy M $. Let $ D_0 $ be a disk representing a special
vertex of $ \D(M,C) $, with $ \bdy D_0 $ in $ \bdy_i $, and let $
\Sigma $ be a sphere splitting off a $ 3 ${\hy}punctured sphere $ P $
containing $ D_0 $ and having $ \bdy_s $ and $ \bdy_i $ as its other
two boundary spheres. As in (b), let $ \DS'(M,C) $ be the subcomplex
of $ \DS(M,C) $ consisting of simplices having no special vertices. We
may deformation retract $ \DS'(M,C) $ onto the link of $ D_0 $ by
putting systems $ S $ representing simplices of $ \DS'(M,C) $ into
normal form with respect to a maximal sphere system containing $
\Sigma $ and then eliminating the intersections of $ S $ with an arc
in $ P $ joining $ \bdy_i $ to $ \bdy_s $ by pushing across $ \bdy_s
$. The intersections of $ S $ with $ P $ can be disks with boundary on
$ \Sigma $, cylinders joining $ \Sigma $ to $ \bdy_i $, or spheres
parallel to $ \bdy_i $.  \figref{fig5} shows which combinations of
these surfaces are possible (along with parallel copies of these
surfaces).

\begin{figure}[ht!]\small\anchor{fig5}
\SetLabels (.015*.2) $\Sigma$\\ (.279*.2) $\Sigma$\\ (.54*.2)
$\Sigma$\\ (.804*.2) $\Sigma$\\ (.108*.69) $\bdy_i$\\ (.108*.26)
$\bdy_s$\\ (.375*.69) $\bdy_i$\\ (.375*.26) $\bdy_s$\\ (.636*.69)
$\bdy_i$\\ (.636*.26) $\bdy_s$\\ (.899*.69) $\bdy_i$\\ (.899*.26)
$\bdy_s$\\ \endSetLabels
\centerline{\AffixLabels{\epsfbox{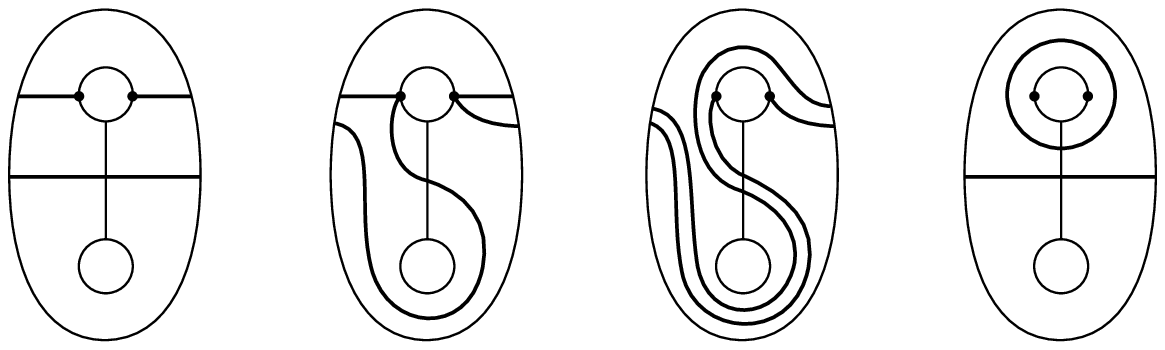}}}
\nocolon\caption{}\label{fig5}
\end{figure}

The rest of the argument then follows the reasoning in (b). The proof
works equally well with $ \D $ in place of $ \DS $.

To deduce (1) from (a), (b), and (c) there are three cases:

\items
\item[(i)] If $ n>0 $ or $ k>1 $ then $ \Sph(M) \ne \emp $, so (b)
gives the reduction to (a) when $ t=0 $ and (c) gives the reduction to
(a) when $ t > 0 $. The conclusion in both cases is that $ \DS(M,C) $
is contractible.

\item[(ii)] Suppose $ n = k = 0 $. If $ t = 0 $ then $ \DS(M,C) =
\Sph(M) \ne \emp $ for $ s \ge $ 4 so (b) gives the reduction to the
case $ s=4 $ and we conclude that $ \DS(M,C) $ is $
(s-5)${\hy}connected. If $ t =1 $ then $ \Sph(M) \ne \emp $ for $ s
\ge 3 $ so (b) gives a reduction to the case $ s = 3 $ and $ \DS(M,C)
$ is $ (s-4)${\hy}connected. If $ t \ge 2 $ then $ \Sph(M) \ne \emp $
so (c) gives a reduction to the case $ s = t $ and we see that $
\DS(M,C) $ is contractible.

\item[(iii)] Suppose $ n = 0 $ and $ k = 1 $. Then $ \Sph(M) \ne \emp
$ for $ s \ge 2 $. Reasoning as in (ii) we see that $ \DS(M,C) $ is
contractible if $ t \ge 2 $ and $ (s-3) ${\hy}connected if $ t \le 1
$. Note that $ s-3 = 2k+s-5 = k+s-4 $ when $ k = 1 $.  \enditems

This finishes the proof of (1). \end{proof}

\begin{proof}[Proof of (2)] Let $ \Sigma $ be a system of spheres that splits off all the $ S^1 \times D^2 $ summands of $M$ whose boundary tori contain circles of $ C $. Doing surgery along this $ \Sigma $ as in the first step of the proof of (1) and discarding all spheres produced by the surgeries gives a deformation retraction of $\D(M,C) $ into the subcomplex consisting of disk systems disjoint from $ \Sigma $. This subcomplex is the join of copies of the disk complex for $ M^1_{0,1} $, which is homeomorphic to the real line. Namely, up to diffeomorphism of $ M^1_{0,1} $ fixing the boundary torus there is a unique nontrivial disk, and applying Dehn twists on the boundary torus gives an infinite sequence of disks, the vertices of the disk complex. There is an edge joining each vertex to each of its two neighbors in the infinite sequence, and there are no simplices of higher dimension, so the disk complex for $ M^1_{0,1} $ is a line. In particular it is contractible, so $ \D(M,C) $ is contractible as well. 

This argument does not cover the special case $ M = M^0_{0,1} $, but
in this case $ \D(M,C) $ is just a point. \end{proof}

\begin{proof}[Proof that (3) implies (4)] For $n=0$ the result holds vacuously, and for $ n=1 $ it is true since $ \D_c(M,C) $ is nonempty if $ n > 0 $ and $ C \ne \emp $. So we may assume $ n \ge 2 $. To show that $ \D_c(M,C) $ is $ (n-2)${\hy}connected,
start with a map $ f \col S^m \rar \D_c(M,C) $, $ m \le n-2 $. By (3)
this extends to a map $ F \col D^{m+1} \rar \D(M,C) $ since $ 2n+k+s-4
\ge n-2 $ if $ n \ge 2 $. We may assume $ F $ is simplicial with
respect to some triangulation of $ D^{m+1} $. Each $p${\hy}simplex $
\sigma$ of $ D^{m+1} $ then maps to a disk system $ D = D_0 \cup
\cdotss \cup D_q $ with $ q \le p $. Choose $ \sigma $ of maximal
dimension $ p $ such that the system $ D $ is purely separating,
meaning that each $ D_i $ separates the complement of the other
disks. This is equivalent to saying that all edges in the dual graph
of $ D $ have distinct endpoints. This implies that $ \sigma $ is not
contained in $ \bdy D^{m+1} $ since $ F $ maps this to systems with
connected complement, having dual graph a wedge of circles. Thus the
link $ L_{\sigma} $ of $ \sigma $ in the triangulation of $ D^{m+1} $
is a sphere $ S^{m-p} $. The maximality of $ \sigma $ implies that any
simplex properly containing $ \sigma $ is sent by $ F $ to a disk
system whose dual graph is obtained from the dual graph of $ D $ by
attaching loops at some vertices. This implies that $ F $ maps $
L_{\sigma} $ to the join
$$
J_{\sigma} = \D_c(M_1,C_1) * \cdotss * \D_c(M_d,C_d)
$$
where the $ M_i $'s are the components of the manifold obtained by
splitting $ M $ along $ D $ and the $ C_i $'s are the preimages of $ C
$ in these components. Each $ (M_i,C_i) $ is smaller than $ (M,C) $
with respect to the lexicographic ordering on the $4${\hy}tuples $
(n,k,s,u) $ where $ u $ denotes the number of circles in $ C $. So by
induction $ \D_c(M_i,C_i) $ is $ (n_i-2)${\hy}connected. Since the
quantity `connectivity-plus-two' is additive for joins by \cite{M}, it
follows that $ J_{\sigma} $ is $ \bigl(\sum_i n_i -
2\bigr)${\hy}connected. We have $ \sum_i n_i \ge n - q \ge n-p $ since
the system $ D $ of $ q+1 $ disks separates $ M $, so splitting $ M $
along $ D $ can decrease the rank of the dual graph by at most $ q
$. Thus $ J_{\sigma}$ is at least $ (n-p-2)${\hy}connected. Since $
L_{\sigma} $ is a sphere of dimension $ m-p $ we deduce that if $ m
\le n-2 $, the map $ F \col L_{\sigma} \rar J_{\sigma} $ extends to a
map $ G \col D^{m-p+1} \rar J_{\sigma} $. The star of $ \sigma $ in $
D^{m+1} $ is the join $ \sigma * L_{\sigma}$ but it can also be
regarded as $ \bdy \sigma * D^{m-p+1} $ for a disk $ D^{m-p+1} \subset
D^{m+1} $ with $ \bdy D^{m-p+1} = L_{\sigma} $. Thus we can redefine $
F $ on the interior of the star of $ \sigma $ by taking the join of
the restriction of $ F $ to $ \bdy \sigma $ with $ G $. For the new $
F $ there are no purely separating simplices of dimension $ p $ or
greater in the star. Hence after a finite number of modifications like
this we eliminate all purely separating simplices, thereby making $ F
$ into a map to $ D_c(M,C) $, finishing the proof that $ D_c(M,C) $ is
$ (n-2) ${\hy}connected. \end{proof}

\begin{proof}[Proof that (4) implies (5)] This is very similar to the preceding proof. The case $ n = 0 $ is easy since $ \D_c(M,C) $ is nonempty, so we may assume $ n \ge 1 $. We start with a map $ S^m \rar \D_c(M,C) $ and extend it to a map $ D^{m+1} \rar \D(M,C) $ using (2). As before, we look at a maximal simplex $ \sigma $ mapping to a purely separating disk system $ D $. The inequality $ \sum_i n_i \ge n-p $ now becomes $ \sum_i n_i \ge n-p +1 $ since splitting along the first disk in a disk system just changes a boundary torus into a sphere. The join $ J_{\sigma} $ is $  \bigl(\sum_i n_i -2 \bigr)${\hy}connected using (4) instead of induction. Then if $ m \le n-1 $ the connectivity of $ J_{\sigma} $ is at least $ \sum_i n_i -2 \ge n - p + 1 - 2 \ge m-p $ so we can finish the proof as before. \end{proof}

\begin{proof}[Proof of (3)] It is natural to try the same approach as in the two previous proofs. Starting with a map $ f \col S^m \rar \D(M,C) $ we can extend it to a map $ F \col D^{m+1} \rar \DS(M,C) $ using (1), and we may assume this is a simplicial map. We wish to modify $ F $ so as to eliminate all sphere systems in its image, so one's first guess would be to take $\sigma$ to be a maximal dimension simplex in $ D^{m+1} $ such that $ F(\sigma) $ is a sphere system. This does not work, however. Here is an example of what can go wrong. Suppose we take $ m=0 $ and we let $ C $ consist of two circles in one sphere of $ \bdy M $. Then there are nontrivial disks $ D_0 $ and $ D_1 $ lying in a neighborhood of this boundary sphere and having boundaries on the two circles of $ C $. These two disks correspond to two distinct vertices of $ \D(M,C) $, and the easiest way to join these two vertices by a path in $ \DS(M,C) $ is to interpolate the vertex formed by a sphere $ S $ parallel to this boundary component of $ M $, so there are edges from $ D_0 $ to $ S $ and from $ S $ to $ D_1 $. Then $\sigma$ would be this vertex $ S $, but there is no way to modify $ F $ on the star of $ \sigma $ by connecting $ D_0 $ and $ D_1 $ by a path of disk systems disjoint from $ S $ since the only such disk systems are $ D_0 $ and $ D_1 $.

So a more subtle approach is needed. Fortunately a model for what to
do in an analogous situation in one lower dimension is available in
Theorem~2.7 of \cite{I}, and we will follow this model closely.

Recall that $ t $ is the number of boundary spheres of $ M $
intersecting $ C $. Note that we can assume that $ t > 0 $ by (2), and
then that $ t = s $ by (1)(c).

For a $p${\hy}simplex $ \sigma $ of $ D^{m+1} $ mapped by $ F $ to a
sphere system $ S $ let $ X = X_{\sigma} $ be the closure of the union
of the components of $ M - S $ that meet $ C $, and let $ Y =
Y_{\sigma} $ be the closure of $ M - X $. Call $ \sigma $ {\em
regular\/} if no spheres of $ S $ lie in the interior of $ Y $. Every
simplex $ \sigma $ of $ D^{m+1} $ mapping to a sphere system contains
a subsimplex which is regular, obtained by discarding vertices mapping
to spheres in the interior of $ Y $. Such vertices are called {\em
superfluous\/}. We will assume from now on that $ \sigma $ is
regular. Let $ Y_1, \cdotss ,Y_c $ be the components of $ Y $ and let
$ X_1, \cdotss, X_d $ be the components of the manifold obtained from
$ X $ by splitting along the spheres of $ S $.

Among the various regular simplices of $ D^{m+1} $ suppose that $
\sigma $ is chosen so that the lexicographically ordered pair $
(c(Y_{\sigma}),p) $ is maximal, where the complexity $ c(Y_{\sigma}) $
is defined as the number of components in the complement of a maximal
sphere system in $ Y_{\sigma} $. (This does not depend on the maximal
sphere system chosen, by the unique prime decomposition theorem for
$3${\hy}manifolds.) With this maximality assumption on $ \sigma $ we
claim that $ F $ takes the link $ L_{\sigma} $ of $ \sigma $ in $
D^{m+1} $ to the join
$$
J_{\sigma} = \D(X_1,C_1) * \cdotss * \D(X_d,C_d)* \Sph(Y_1)* \cdotss *
\Sph(Y_c)
$$
where $ C_i $ is the part of $ C $ in $ X_i $. Indeed, let $ v $ be a
vertex in $ L_{\sigma} $, so $ F(v) $ is a disk or sphere in the
complement of the system $ S = F(\sigma) $. If $ F(\sigma) $ lies in
some $ Y_j $ it must be a sphere defining a vertex of $ \Sph(Y_j) $
and hence in $ J_{\sigma} $. The other possibility is that $ F(\sigma)
$ lies in some $ X_i $, where it could be a disk or a sphere, and we
need to show it must be a disk. So suppose it is a sphere. If it does
not split off a piece of $ X_{\sigma} $ disjoint from $ C $ then we
can enlarge $ \sigma $ to include $ v $ without changing $ X_{\sigma}
$ or $ Y_{\sigma} $, contradicting the maximality of $
(c(Y_{\sigma}),p) $.  On the other hand, if $ F(\sigma) $ does split
off a piece of $ X_{\sigma} $ disjoint from $ C $ then by enlarging $
\sigma $ to include $ v $ and then discarding any resulting
superfluous vertices we obtain a new $ \sigma' $ with $ Y_{\sigma'} $
the union of $ Y_{\sigma} $ and the piece split off from $ X_{\sigma}
$, so $c(Y_{\sigma'}) > c(Y_{\sigma}) $, contradicting maximality
again.

So we have $ F \col L_{\sigma} \rar J_{\sigma} $. We calculate now the
connectivity $ J_{\sigma} $. If any $ Y_j $ has $ n(Y_j) >0 $ or $
k(Y_j) > 1 $ then $ \Sph(Y_j) $ is contractible by (1), hence $
J_{\sigma} $ is also contractible. So we may assume $ n(Y_j) = 0 $ and
$ k(Y_j) \le 1 $ for each $ j $. For the connectivity of the join $
J_{\sigma} $ we have
$$
\conn(J_{\sigma}) = {\textstyle\sum_i} \conn\bigl(\D(X_i,C_i)\bigr) +
{\textstyle\sum_j} \conn\bigl(\Sph(Y_j)\bigr) + 2c + 2d -2
$$
By induction on the lexicographically ordered triple $ (n,k,s) $ we
know the connectivities in the first summation, and (1) gives the
connectivities in the second summation, so if the number of $ Y_j$'s
with $ k(Y_j) = 1 $ is $ a $, we get
\begin{align*}
\conn(J_{\sigma}) = 2{\textstyle\sum_i} n(X_i) + &{\textstyle\sum_i}
 k(X_i) + {\textstyle\sum_i} s(X_i) -4d\\ &+ {\textstyle\sum_j} s(Y_j)
 - 5c + 2a + 2c + 2d - 2 \end{align*} We can start the induction with
 $ (n,k,s) = (0,0,3) $ when $ \D(M,C) $ is nonempty if $ C \ne \emp
 $. If $ F(\sigma) $ is a $q${\hy}simplex, consisting of $ q+1 $
 spheres, then by splitting $ M $ along these spheres and counting the
 number of boundary spheres in the resulting manifold we see that $
 s(M) + 2q + 2 = \sum_i s(X_i) + \sum_j s(Y_j) $. Looking at the dual
 graph to the system $ F(\sigma) $ we see that $ n(M) = \sum_i n(X_i)
 + q + 1 - (c+d-1) $ since it takes $ c+d-1 $ edges to connect
 together the $ c + d $ subgraphs corresponding to the $ X_i $'s and $
 Y_j $'s and we assumed $ n(Y_j) = 0 $. Also we have $ k(M) = \sum_i
 k(X_i) + a $. Substituting into the preceding displayed formula and
 letting $ n = n(M) $, $ k = k(M) $, $ s = s(M) $, we see that the
 connectivity of $ J_{\sigma} $ is
\begin{align*}
 (2n - & 2q -4 + 2c + 2d) + (k-a) - 4d + (s + 2q+2) - 5c + 2a + 2c +
2d -2 \cr &= 2n + k + s -4 + a - c \end{align*}

Note that $ q + 1 \ge 3(c - a) $ since each $ Y_j$ with $ k(Y_j)= 0 $
must have at least three boundary spheres and we assumed that no
boundary sphere of $ M $ was disjoint from $ C $. This implies that $
q \ge c - a $. If we assume that $ m \le 2n + k + s -4 $, then
\begin{align*}
\dim L_{\sigma} = m-p \le m-q &\le 2n + k + s -4 - q \cr &\le 2n + k +
s -4 + a - c = \conn(J_{\sigma}) \end{align*}

Hence if $ m \le 2n+k+s-4 $ we can modify $ F $ in the star of $
\sigma $ as before, taking the join of its values on $\bdy\sigma$ with
an extension of $ F \col L_{\sigma} \rar J_{\sigma} $ to a map $
D^{m-p+1} \rar J_{\sigma} $.

It remains to check that this modification improves the situation,
namely that the maximum value of $ (c(Y_{\sigma}),\dim\sigma) $ has
not increased, while the number of simplices $ \sigma $ realizing this
maximum value has actually decreased. We only need to check simplices
$ \tau $ in $ D^{m-p+1} * \bdy\sigma $ that are sent by the new $ F' $
to sphere systems, and we need only consider such $ \tau $ that are
regular. The simplex $ \tau $ is the join of a simplex $ \alpha $ in $
D^{m-p+1} $ and a simplex $ \beta $ in $ \bdy\sigma $, where $ \alpha
$ or $ \beta $ is taken to be empty when $ \tau $ lies in $\bdy\sigma
$ or $ D^{m-p+1} $. Since $ F'(\tau) $ is a sphere system, so are $
F'(\alpha) $ and $ F'(\beta) $. By construction $ F'(\alpha) $ is a
simplex in $ J_{\sigma} $, so since it is a sphere system it must lie
in $ \Sph(Y_1) * \cdotss * \Sph(Y_c) $ and hence it must consist of
spheres in $ Y_{\sigma} $. On the other hand, $ F'(\beta) $ is a
subsimplex of $ F(\sigma) $.  These two facts imply that $ Y_{\tau}
\subset Y_{\sigma} $. Hence $ c(Y_{\tau}) \le c(Y_{\sigma})
$. Equality can occur only if $ Y_{\tau} = Y_{\sigma} $, in which case
$ \alpha = \emp $ since $ \tau $ is regular. Then the simplex $ \tau =
\beta $ is a proper face of $ \sigma $ and hence has smaller
dimension. So in all cases $ (c(Y_{\tau}),\dim\tau) <
(c(Y_{\sigma}),\dim\sigma) $. \end{proof}

This finishes the proof of Theorem \ref{proc3.1}. We will also need a
similar result with spheres in place of disks. Let $\Sph_c(M)$ be the
subcomplex of $\Sph(M)$ consisting of sphere systems with connected
complement and let $\Sph^\pm_c(M)$ be the complex of oriented sphere
systems with connected complement. The simplices of $ \Sph^\pm_c(M) $
thus correspond to systems of spheres each of which has a chosen
orientation, or equivalently, a chosen normal orientation. (There is
no need to consider oriented disk complexes since disks can be
canonically oriented by choosing orientations for the circles of $ C
$.)

\begin{prop}\label{proc3.2}
$\Sph_c(M)$ and $\Sph^\pm_c(M)$ are $(n-2)${\hy}connected.\end{prop}

This was shown in Proposition 3.1 of \cite{H3} in the case $ k = 0
$. The generalization to $ k > 0 $ is straightforward, and we give the
proof here mainly to correct a small error at the beginning of the
proof in \cite{H3}.

\begin{proof}  For the case of $\Sph_c(M)$ this will follow from (1) by the same argument used to show that (3) implies
(4).  The result is trivially true for $n=0$, so we may assume that
$n>0$.  Consider a map $f\col S^m\rar \Sph_c(M)$ with $m\le
n-2$. Since $ \Sph(M) $ is contractible, this extends to a map $F\col
D^{m+1}\rar \Sph(M)$. Let $\sigma$ be a purely separating simplex in
$D^{m+1}$ of maximal dimension $p$ with image the sphere system $S =
S_0\cup\cdotss\cup S_q$ where $q\le p$. Then $F$ maps the link
$L_\sigma = S^{m-p}$ to $J_\sigma = \Sph_c(M_1)*\cdotss * \Sph_c(M_d)$
where the $M_i$'s are the manifolds obtained by cutting $M$ along
$S$. By induction on the triple $(n,k,s)$, the connectivity of
$\Sph_c(M_i)$ is $n_i-2$ and thus $J_\sigma$ is $(\sum_i n_i
-2)${\hy}connected. Now $\sum_i n_i \ge n-q$ as the system $S $
separates $ M $. So $\sum_i n_i -2 \ge n-q-2 \ge m-p$ as $m\le n-2$
and $q\le p$. Hence $F\col L_\sigma= S^{m-p} \rar J_\sigma$ can be
extended to the disc $D^{m-p+1}$. We use this to modify $F$ on the
interior of the star of $\sigma$.  After a finite number of such
modifications, $F$ becomes a map with image in $\Sph_c(M)$, so $
\Sph_c(M) $ is $ (n-2)${\hy}connected.

We want to deduce that $\Sph_c^\pm(M)$ is also $(n-2)${\hy}connected.
For this, choose a positive orientation on the spheres of $\Sph_c(M)$
so that $\Sph_c(M)$ is identified with the subcomplex $\Sph_c^+(M)$ of
$\Sph^\pm_c(M)$ of spheres with positive orientation. It will suffice
to take a map $f\col S^m\rar \Sph_c^\pm(M)$ with $m\le n-2$ and
homotope this to have image in $ \Sph^+_c(M) $. Let $\sigma$ be a
simplex of $S^m$ of maximal dimension $p$ whose image under $ f $ is a
system $ S = S^-_0\cup\cdotss\cup S^-_q$ consisting entirely of
negatively oriented spheres. The link $L_\sigma=S^{m-p-1}$ then maps
to $\Sph_c^+(M')$, where $M'$ is $M$ cut along $S$. The manifold $M'$
has $n-q-1$ copies of $S^1\times S^2$ in its decomposition as $M - S$
is connected. So $S^+_c(M')$ is $(n-q-3)${\hy}connected. Now $m-p-1\le
n-q-3$ as $m\le n-2$ and $q\le p$. Hence we can modify $f$ on the
interior of the star of $\sigma$ as before. Note that the new $ f $ is
homotopic to the old one by coning. After a finite number of such
modifications we obtain an $ f $ having image in $ \Sph^+_c(M) $.
\end{proof}

We conclude this section with two results that will be needed for the
spectral sequence argument. We use the notation $ A(M) $ for the
quotient of $ \pi_0\Dif(M) $ by twists along $ 2 ${\hy}spheres. This
quotient acts on all the disk and sphere complexes for $ M $ since
twists along $2${\hy}spheres act trivially, as they obviously preserve
the homotopy classes of embedded spheres and disks.

\begin{prop}\label{proc3.3}
If $ C $ consists of a single circle, then $ A(M) $ acts transitively
on the $ p ${\hy}simplices of $ \D_c(M,C) $ for each $ p $. The same
is true also for $ \Sph_c(M) $ and $ \Sph^\pm_c(M) $.
\end{prop}

\begin{proof} First consider the case that $ C $ is contained in a torus boundary
component. Let $ D = D_0 \cup \cdotss \cup D_p $ be a disk system
representing a $ p ${\hy}simplex of $ \D_c(M,C) $.  Thinking of $ C $
as a longitude of the torus $ T $ of $ \bdy M $ that it lies in, take
a meridian circle of $ T $ and push this slightly into the interior of
$ M $. Orienting this meridian circle then determines an ordering of
the disks $ D_i $, which we may assume agrees with the ordering by
increasing subscripts. If $ D' = D'_0 \cup \cdotss \cup D'_p $ is
another such disk system, we want to find a diffeomorphism $ M \rar M
$ fixed on $ \bdy M $ and taking each $ D_i $ to $ D'_i $.

Let $ N $ be the manifold obtained from $ M $ by splitting along $ D
$. Splitting along $ D_0 $ replaces the $ S^1 \times D^2 $ summand of
$ M $ containing $ T $ by a $ D^3 $ summand, and then, since $ M - D $
is connected, splitting along each subsequent $ D_i $ replaces an $
S^1 \times S^2 $ summand by another $ D^3 $ summand. The same is true
when we split along $ D' $ to produce another manifold $ N' $. By the
uniqueness of prime decompositions for compact orientable
$3${\hy}manifolds, $ N $ and $ N' $ are diffeomorphic, and the
diffeomorphism can be chosen to be the identity on the components of $
\bdy M $ other than $ T $. The diffeomorphism can also be chosen to
permute the new boundary spheres in any way we like, so it can be
chosen in such a way that it passes down to a quotient diffeomorphism
$ M \rar M $ taking each $ D_i $ to $ D'_i $ and fixed on $ \bdy M $.

The case that $ C $ is contained in a sphere boundary component is
covered by the same argument since in this case splitting along disks
again replaces $ S^1 \times S^2 $ summands by $ D^3 $ summands. For $
\Sph_c(M) $ and $ \Sph^\pm_c(M) $ splitting along spheres replaces $
S^1 \times S^2 $ summands by pairs of $ D^3 $ summands, so a similar
argument applies. \end{proof}

For a system $ S $ of disks and spheres in $ M $ let $ A(M,S) $ be the
subgroup of $ A(M) $ represented by diffeomorphisms restricting to the
identity on $ S $. If $ M' $ denotes $ M $ split along $ S $, there is
thus a surjection $ A(M') \rar A(M,S) $.

\begin{prop}\label{proc3.4}
The surjection $ A(M') \rar A(M,S) $ is an isomorphism.\end{prop}

\begin{proof}  Choose a system $ \Sigma $ in $ M' $ consisting of spheres 
that split off all prime summands together with spheres in $ S^1
\times S^2 $ summands and disks in $ S^1 \times D^2 $
summands. Represent an element of the kernel of $ A(M') \rar A(M,S) $
by a diffeomorphism $ f \in \Dif(M) $ that is isotopic to the identity
and restricts to the identity on $ S $. The system $ f(\Sigma) $ is
then isotopic to $ \Sigma $ in $ M $, hence in $ M' $ as shown in the
proof of Theorem~\ref{proc2.2}. So we may assume $ f $ takes each
component of $ \Sigma $ to itself. We may also assume $ f $ is the
identity on disk components of $ \Sigma $ since it is the identity on
their boundaries. On the spheres in $ \Sigma $ that are separating, $
f $ must preserve orientation since it cannot switch the sides of
these spheres, and $ f $ must also preserve orientations of the
nonseparating spheres in $ \Sigma $ since it induces the identity on
homology. Thus we may assume $ f $ is the identity on $ \Sigma $ as
well as on $ S $. Modulo twists along $2${\hy}spheres we can then
deform $ f $, staying fixed on $ \Sigma $, to be the identity on all
of $ M $ except the prime summands that are neither $ S^1 \times S^2 $
nor $ S^1 \times D^2 $, since an orientation-preserving diffeomorphism
of a punctured $ S^3 $ that takes each puncture to itself is isotopic
to the identity. For the remaining prime summands we appeal to the
theorem in Section 5 of \cite{HL} which has as a corollary the fact
that the product of the mapping class groups of these prime summands
injects into the mapping class group of $ M $. \end{proof}

\section{Proof of stability}

We make a small change of notation now and let $ M^s_{n,k} $ be the
connected sum of $ n $ copies of $ S^1 \times S^2 $, $ k $ copies of $
S^1 \times D^2 $, $ s $ balls, and perhaps also another manifold $ N $
that is a sum of prime factors satisfying the conditions imposed at
the beginning of Section 2. Then we have the following stabilization
result for the groups $ A^s_{n,k} $, which are the mapping class
groups of the manifolds $ M^s_{n,k} $ with twists along $ 2
${\hy}spheres factored out:

\begin{thm}\label{proc4.1}
For all $ k \ge 0 $ and $ s \ge 1 $ the natural stabilization maps
\begin{align*}
& \hbox {\rm (1)\enskip} \ H_i(A^s_{n,k}) \rar H_i(A^{s+1}_{n,k})
\qquad\qquad\qquad\qquad\qquad\qquad\qquad\cr & \hbox {\rm (2)\enskip}
\ H_i(A^s_{n,k}) \rar H_i(A^s_{n,k+1}) \cr & \hbox {\rm (3)\enskip} \
H_i(A^s_{n,k}) \rar H_i(A^s_{n+1,k}) \end{align*}

are isomorphisms whenever $ n \ge 3i + 3 $.\end{thm}

 These isomorphisms will follow from three other similar
isomorphisms. Consider the following three maps, for $ s \ge 0 $:

\items
\item[\lowsqdot] $\alpha\col A_{n,k}^{s+2}\rar A_{n+1,k}^{s+1}$,
induced by the map $ M_{n,k}^{s+2}\rar M_{n+1,k}^{s+1}$ identifying
disks in each of the last two boundary spheres of $M$, or
equivalently, attaching a $1${\hy}handle $ D^1 \times D^2 $ joining
these two boundary spheres.
\item[\sqdot] $\beta\col A_{n,k}^{s+2}\rar A_{n+1,k}^s$, induced by
the map $M_{n,k}^{s+2}\rar M_{n+1,k}^s$ identifying the two last
boundary spheres of $M$, or equivalently, joining them by a product $
S^2 \times I $.
\item[\lowsqdot] $ \gamma\col A^{s+1}_{n,k} \rar A^s_{n,k+1} $,
induced by the map $ M^{s+1}_{n,k} \rar M^s_{n,k+1} $ identifying two
disks in the last boundary sphere of $ M^{s+1}_{n,k} $, or
equivalently, attaching a $1${\hy}handle $ D^1 \times D^2 $ from this
boundary sphere to itself.  \enditems

What we will actually prove is that for all $ k \ge 0 $ and $ s \ge 0
$ the induced maps
\begin{align*}
& \hbox {\rm (A)\enskip} \ \alpha_i\col H_i(A_{n,k}^{s+2}) \rar
H_i(A_{n+1,k}^{s+1}) \qquad\qquad\qquad\qquad\qquad\cr & \hbox {\rm
(B)\enskip} \ \beta_i\col H_i(A_{n,k}^{s+2}) \rar H_i(A_{n+1,k}^{s})
\cr & \hbox {\rm (C)\enskip} \ \gamma_i\col H_i(A_{n,k}^{s+1}) \rar
H_i(A_{n,k+1}^{s}) \end{align*}

are surjective when $n\ge 3i$ and isomorphisms when $n\ge 3i+2$. To
see how this implies the theorem, we use two other maps:

\items
\item[\sqdot] $ \delta \col A^{s+1}_{n,k} \rar A^s_{n,k} $, induced by
filling in the last boundary sphere with a ball.

\item[\lowsqdot] $ \varepsilon \col A^s_{n,k} \rar A^{s+1}_{n,k} $,
induced by attaching a $2${\hy}handle $ D^2 \times D^1 $ to the last
boundary sphere, if $ s \ge 1 $.  \enditems

Since $ \beta $ is the composition $ \delta\alpha $, we deduce that $
\delta_i\col H_i(A^{s+1}_{n,k}) \rar H_i(A^s_{n,k}) $ is an
isomorphism for $ n \ge 3i+3 $. The composition $ \delta\varepsilon $
is the identity, so we conclude that $ \varepsilon $ induces an
isomorphism on homology:

\items
\item[(1)] The stabilization $ H_i(A^s_{n,k}) \rar H_i(A^{s+1}_{n,k})
$ is an isomorphism for $ n \ge 3i+3 $ and $ s \ge 1 $.  \enditems

The composition $ \gamma\varepsilon $ comes from attaching disjoint $
1 ${\hy} and $ 2 ${\hy}handles to the last boundary sphere, so this is
the standard stabilization $ A^s_{n,k} \rar A^s_{n,k+1} $. Thus we
have:

\items
\item[(2)] The stabilization $ H_i(A^s_{n,k}) \rar H_i(A^s_{n,k+1}) $
is an isomorphism for $ n \ge 3i+3 $ and $ s \ge 1 $.  \enditems

Finally, the standard stabilization $ A^s_{n,k} \rar A^s_{n+1,k} $ is
$ \beta\varepsilon^2 $, so:

\items
\item[(3)] The stabilization $ H_i(A^s_{n,k}) \rar H_i(A^s_{n+1,k}) $
is an isomorphism for $ n \ge 3i+3 $ and $ s \ge 1 $.  \enditems

\begin{cor}\label{proc4.2}
The quotient map $\Aut(F_n) \rar \Out(F_n)$ induces an isomorphism on
$ H_i $ for $ n \ge 3i+3 $.\end{cor}

\begin{proof} The map $\Aut(F_n) \rar \Out(F_n)$ is a special case of $ \delta $ as $ A^1_{n,0} = \Aut(F_n)
$ and $ A^0_{n,0} = \Out(F_n) $, and we have just seen that $ \delta $
induces an isomorphism on $ H_i $ in the stated range. \end{proof}

\begin{proof}[Case (A)] We do an induction on $i$, using the relative spectral sequence
argument of \cite{V}.  The statement is clearly true for $i=0$. The
induction hypothesis will be that the following two statements are
true for all $ q < i $:

\items
\item[$(a_q)$] $\alpha_q\col H_q(A_{n,k}^{s+2}) \rar
H_q(A_{n+1,k}^{s+1})$ is surjective for all $n\ge 3q$.

\item[$(b_q)$] $\alpha_q\col H_q(A_{n,k}^{s+2}) \rar
H_q(A_{n+1,k}^{s+1})$ is an isomorphism for all $n\ge 3q+2$.
\enditems

 We want to show that this implies $(a_i)$ and $(b_i)$.

To ease the notation we write the stabilization $ A^{s+2}_{n,k} \rar
A^{s+1}_{n+1,k} $ as $ G_n \rar G_{n+1} $. Consider the action of
$G_{n+1} $ on the complex $ X_{n+1} = \D_c(M_{n+1,k}^{s+1},C) $, where
$C$ is a single circle in the last sphere component of $ \bdy M $. The
stabilizer of a vertex is $ G_n $, and more generally the stabilizer
of a $ p ${\hy}simplex is $ G_{n-p} = A^{s+p+2}_{n-p,k} $.  If we take
$ X_n $ to be $ \D_c(M_{n,k}^{s+2},C) $ with $ C $ again a single
circle in the last boundary sphere, then the vertex stabilizer $ G_n $
acts on $ X_n $, and the map $ (M_{n,k}^{s+2},C) \rar
(M_{n+1,k}^{s+1},C) $ that induces $ \alpha $ also induces an
embedding of $ X_n $ in $ X_{n+1} $ preserving the actions of $ G_n $
and $ G_{n+1} $. Thus if $(C_*,\del_*)$ and $(C'_*,\del'_*)$ are the
augmented chain complexes of $X_n$ and $X_{n+1}$ respectively, then we
have a map of double complexes
$$
i\col E_qG_n\otimes_{G_n} C_p \Rar E_q G_{n+1}\otimes_{G_{n+1}} C'_p
$$ 
where $(E_*G_n,d_*)$ is a free $G_n${\hy}resolution of $\Z$ and
similarly for $ (E_*G_{n+1},d_*)$.

With mapping cones in mind, consider now the double complex
$$
(E_{q-1}G_n\otimes_{G_n} C_p) \oplus (E_q G_{n+1}\otimes_{G_{n+1}}
C'_p)
$$
with vertical boundary maps $ {(a \smalltensor b, a' \smalltensor b')
\mapsto \bigl(-da \smalltensor b, d'a' \smalltensor b' + i(a
\smalltensor b)\bigr)} $ and horizontal boundary maps $ (a
\smalltensor b, a' \smalltensor b') \mapsto (-1)^q(a \smalltensor \bdy
b,a' \smalltensor \bdy'b') $.  There are two spectral sequences
associated to this double complex.  The first spectral sequence,
arising from the horizontal filtration, converges to 0 for $p\le n-2$
since $H_p(C_*)=0$ for $p\le n-2$ and $H_p(C'_*)=0$ for $p\le n-1$ by
Theorem \ref{proc3.1}. So in the second spectral sequence, associated
to the vertical filtration, we have $E^\infty_{p,q}=0$ for $p+q\le
n-2$.  The columns in the double complex are the mapping cones of the
chain maps $i\col E_*G_n\otimes_{G_n} C_p \rar
E_*G_{n+1}\otimes_{G_{n+1}} C'_p$. So the $E^1$ term of the second
spectral sequence is
$$
E^1_{p,q}=H_q(E_*G_{n+1}\otimes_{G_{n+1}} C'_p,E_*G_n\otimes_{G_n}
C_p)
$$
By Proposition \ref{proc3.3} the actions of $G_n$ and $ G_{n+1} $ on
$X_n$ and $ X_{n+1} $ are transitive on the sets of
$p${\hy}simplices. For the map $X_n \rar X_{n+1}$, each orbit of a
$p${\hy}simplex in $ X_{n+1} $ comes from the orbit of a
$p${\hy}simplex in $ X_n $ if $ p \le \dim X_n = n-1 $. The
stabilizers of $ p ${\hy}simplices in $ X_n $ and $ X_{n+1} $ are $
G_{n-p-1} $ and $ G_{n-p} $ respectively, so Shapiro's Lemma (or its
proof) gives an isomorphism
$$
E^1_{p,q} \cong H_q(G_{n-p},G_{n-p-1}) \quad {\rm for }\ p \le n-1
$$
where the coefficients are untwisted $ \Z $'s since the stabilizers do
not permute the vertices of a simplex.

To prove $(a_i)$, we will show that $H_i(G_{n+1},G_n)=0$ when $n\ge
3i$.  Consider the differential
$$
d^1\col E^1_{0,i}=H_i(G_n,G_{n-1}) \Rar E_{-1,i}^1=H_i(G_{n+1},G_n)
$$
By assumption, $H_q(G_{n-p},G_{n-p-1})=0$ if $q<i$ and $p+q=i$ as
$n-p-1\ge 3i-p-1\ge 3i-3p=3q$ since $p\ge 1$.  As $n\ge 3i$ also
implies $i-1\le n-2$ when $i>0$, the term $E_{-1,i}^1$ must be killed
by differentials and hence the differential $d^1$ above is surjective
as the other differentials originate in trivial groups. Now a diagram
chase in
\vskip-12pt
$$ \SelectTips{cm}{} \xymatrix{H_i(G_n) \ar[r] \ar[d] &
H_i(G_n,G_{n-1}) \ar[d]^{d^1} \ar[r] & H_{i-1}(G_{n-1}) \ar[r]^-\cong
& H_{i-1}(G_n)\\ H_i(G_{n+1}) \ar[r] & H_i(G_{n+1},G_n) & & }$$ shows
that $H_i(G_{n+1},G_n)=0$.  Indeed, given an element $x\in
H_i(G_{n+1},G_n)$ there is an $x'\in H_i(G_n,G_{n-1})$ mapping to
$x$. The top right map is an isomorphism by $(b_{i-1})$ as $n-1\ge
3i-1= 3(i-1)+2$. So $x'$ maps to $0$ in $ H_{i-1}(G_{n-1})$, and thus
there is an $x''\in H_i(G_n)$ mapping to $x'$. Following $ x'' $
around the other two sides of the square gives $ 0 $ since these
groups are part of an exact sequence, so $x=0$. This proves $(a_i)$.

The proof of $(b_i)$ is similar. From the preceding paragraph and the
assumption $n\ge 3i+2$ in $(b_i)$ we have
$E^1_{p,q}=H_q(G_{n-p},G_{n-p-1})=0$ for all $q\le i$ and $p+q=i+1$ as
$n-p-1\ge 3i+2-p-1\ge 3i-3p+3=3q$ since $ p \ge 1 $. From the spectral
sequence we conclude by the same reasoning as before that the
differential
$$
d^1\col E^1_{0,i+1}=H_{i+1}(G_n,G_{n-1}) \Rar
E_{-1,i+1}^1=H_{i+1}(G_{n+1},G_n)
$$
is surjective as $i\le n-2$ when $n\ge 3i+2$.  A diagram chase in
\vskip-10pt
$$ \SelectTips{cm}{} \xymatrix{H_{i+1}(G_n,G_{n-1}) \ar[r]
\ar[d]^{d^1} & H_i(G_{n-1}) \ar[d] & \\ H_{i+1}(G_{n+1},G_n) \ar[r] &
H_i(G_n) \ar[r]^-{\alpha_i} & H_i(G_{n+1})}$$ now shows that
$\alpha_i\col H_i(G_n) \rar H_i(G_{n+1})$ is injective. Indeed, let
$x\in H_i(G_n)$ be in the kernel of $\alpha_i$. Then there exists
$x'\in H_{i+1}(G_{n+1},G_n)$ mapping to $x$. As the left vertical map
is surjective, there is a lift $x''\in H_{i+1}(G_n,G_{n-1})$ of
$x'$. Following $ x'' $ around the other two sides of the square gives
$ 0 $, so $ x = 0 $. This proves $(b_i)$ since $(a_i)$ gives the
surjectivity of~$ \alpha_i $.
\end{proof}

\begin{proof}[Case (B)] The result in this case is shown in \cite{HV3} by a slightly different spectral
sequence argument, and we include a proof here mainly for the reader's
convenience.

As in the previous case, we do an induction on $i$. The statement is
clearly true for $i=0$. Suppose inductively that the following are
valid for all $ q < i $:

\items
\item[$(a_q)$] $\beta_q\col H_q(A^{s+2}_{n,k}) \rar H_q(A^s_{n+1,k}) $
is surjective for all $n\ge 3q$.

\item[$(b_q)$] $\beta_q\col H_q(A^{s+2}_{n,k}) \rar H_q(A^s_{n+1,k})$
is an isomorphism for all $n\ge 3q+2$.  \enditems

We write the stabilization $ A^{s+2}_{n,k} \rar A^s_{n+1,k} $ as $ G_n
\rar G_{n+1} $ and consider the action of $ G_{n+1} $ on $ X_{n+1} =
\Sph^\pm_c(M^s_{n+1,k}) $. A vertex stabilizer is $ G_n $, but
unfortunately the stabilizer of a $ p ${\hy}simplex for $ p > 0 $ is
larger than $ G_{n-p} = A^{s + 2p+2}_{n-p,k} $ since the stabilizer
includes diffeomorphisms that permute the vertices of the
$p${\hy}simplex. To fix this problem one could attempt to refine the
definition of $ \Sph^\pm_c(M^s_{n+1,k}) $ so as to include the data of
an ordering of the spheres in a sphere system. However, it is not
evident that this complex would be highly connected. So instead we do
something slightly different, letting $ X_{n+1} $ be the complex whose
$ p ${\hy}simplices are all the simplicial maps $ \Delta^p \rar
\Sph^\pm_c(M^s_{n+1,k}) $. These maps take vertices to vertices, but
different vertices of $ \Delta^p $ can have the same image, so we are
allowing `degenerate' simplices. (This means that $ X_{n+1} $ is no
longer a simplicial complex, but it is a CW complex with the added
structure of a $\Delta${\hy}complex, in the terminology of \cite{H4}.)
It is a standard fact in algebraic topology that enlarging a
simplicial complex by adjoining degenerate ordered simplices in this
way does not change its homology groups; see for example Theorem~4.6.8
of \cite{S}. The group $ G_{n+1} $ still acts on the new $ X_{n+1} $,
and the action is transitive on vertices with $ G_n $ as vertex
stabilizer, but the action is no longer transitive on
higher-dimensional simplices. The stabilizers of $ p ${\hy}simplices
are the groups $ G_{n-r} $ for $ r \le p $, and this will be good
enough for the proof in case (A) to go through without significant
change.

The action of $ G_{n+1} $ on $ X_{n+1} $ restricts to an action of $
G_n $ on $ X_n $, so as in case (A) we look at the double complex
$(E_{q-1}G_n\otimes_{G_n} C_p) \oplus (E_qG_{n+1}\otimes_{G_{n+1}}
C'_p) $ with boundary maps as before. The spectral sequence coming
from the horizontal filtration converges to 0 when $p\le n-2$, so the
second spectral sequence has $E^\infty_{p,q}=0$ for all $p+q\le n-2$.

The inclusion $ X_n \inc X_{n+1} $ induces a bijection on orbits of $
p ${\hy}simplices for $ p \le n-1 $, so Shapiro's Lemma says that the
term $E^1_{p,q}$ of the spectral sequence obtained from the vertical
filtration is a direct sum of terms $ H_q(G_{n-r},G_{n-r-1}) $ for $ r
\le p $ if $ p \le n-1 $. Again the coefficients are untwisted $ \Z
$'s since stabilizers do not permute vertices of a simplex.

The rest of the argument now proceeds as in case (A) as $n-r-1\ge 3q$
for all $r\le p$ if and only if this holds when $r=p$. \end{proof}

\begin{proof}[Case (C)] In order to keep the same notation as in the previous cases we take the
liberty of writing $ A^{s+1}_{n,k} \rar A^s_{n,k+1} $ as $ G_n \rar
G_{n+1} $.  This is the map induced by the inclusion
$M_{n,k}^{s+1}\inc M_{n,k+1}^s$. We consider the action of $ G_{n+1} $
on $ X_{n+1} = \D_c(M^s_{n,k+1},C) $ where $ C $ is a single circle in
the last boundary torus. Then $ G_n $ is a vertex stabilizer, and this
acts on $ X_n = \D_c(M^{s+1}_{n,k},C) $ where $ C $ is a single circle
in the last boundary sphere, as in case (A). The connectivities of $
X_n $ and $ X_{n+1} $ are $n-2$ and $n-1$ respectively and the
inclusion $M_{n,k}^{s+1}\inc M_{n,k+1}^s$ induces an equivariant map
$X_n\rar X_{n+1}$. The actions are transitive on simplices in each
dimension.  The stabilizer of a $ p ${\hy}simplex for the action of $
G_{n+1} $ on $ X_{n+1} $ is $ G_{n-p} = A^{s+p+1}_{n-p,k} $ and for
the action of $G_n$ on $X_n$ it is
$A_{n-p-1,k}^{s+p+2}=G_{n-p-1}$. Thus the stabilizers are the same as
in (A). More precisely, taking $\tild s=s+1$ in (A), the
$E^1${\hy}pages are the same when $0\le p\le n-1$, so we already know
that the appropriate terms are $0$. The result thus follows by the
same argument as for (A), except that we do not need a full induction
argument, but can instead appeal to the results proved in
(A). \end{proof}

It is worth noting that $ \gamma $ automatically induces injections on
all homology groups since if we follow it by the map coming from
attaching a $ 2 ${\hy}handle to the last boundary torus to surger it
to a sphere, the composition is the identity. One can also see from
the proof of case (C) that an improvement in the stable range for
stabilization with respect to $ n $ and $ s $ would yield an
improvement for stabilization with respect to $ k $.

\end{document}